\documentclass[pdflatex]{article}

\usepackage{hyperref}
\usepackage{authblk}
\usepackage{graphicx}%
\usepackage{multirow}%
\usepackage{amsmath,amssymb,amsfonts}%
\usepackage{amsthm}%
\usepackage[title]{appendix}%
\usepackage[dvipsnames]{xcolor}%
\usepackage{textcomp}%
\usepackage{manyfoot}%
\usepackage{booktabs}%
\usepackage{algorithm}%
\usepackage{algorithmicx}%
\usepackage{algpseudocode}%
\usepackage{listings}%
\usepackage[shortlabels]{enumitem}
\usepackage[a4paper,left=3.25cm,right=3.25cm]{geometry}
\usepackage{xfrac} 
\usepackage{tikz}
\usepackage{tikz-3dplot}
\usetikzlibrary{arrows.meta, calc, backgrounds, positioning, shadows}
\usetikzlibrary{spy}
\usepackage{xinttools} 
\usepackage[caption=false]{subfig}
\usepackage{comment}
\usepackage{textgreek} 

\usepackage{fancybox}
\usepackage{mathtools}
\usepackage{blkarray}
\usepackage{framed}

\usepackage{mathbbol}
\usepackage{lmodern}
\usepackage[normalem]{ulem}

\newtheorem{Th}{Theorem}
\newtheorem{Prop}[Th]{Proposition}

\theoremstyle{definition}

\newcommand{\RR}{\mathbf{R}}

\newcommand{\conv}[2][n]{\underset{#1\rightarrow #2}{\longrightarrow}}
\newcommand{\eq}[2][n]{\underset{#1\rightarrow #2}{\sim}}

\newcommand{\EEE}[1]{\operatorname{\mathbb{E}}\left[\,#1\,\right]}
\newcommand{\VVV}[1]{\operatorname{\mathbb{V}}\left[\,#1\,\right]}

\newcommand{\PPP}[1]{\operatorname{\mathbb{P}}\left(\,#1\,\right)}

\newcommand{\ind}[1]{\mathbb{1}_{#1}\,}

\newcommand{\eps}{\varepsilon}
\newcommand{\WP}{W}

\newcommand{\WPe}{\WP'}

\newcommand{\WPP}{[-\frac{1}{2}, \frac{1}{2}]^2}
\newcommand{\WPPe}{[-\frac{1}{2} - \eps, \frac{1}{2} + \eps]^2}

\def\biglen{20cm} 
\tikzset{
  half plane/.style={ to path={
       ($(\tikztostart)!.5!(\tikztotarget)!#1!(\tikztotarget)!\biglen!90:(\tikztotarget)$)
    -- ($(\tikztostart)!.5!(\tikztotarget)!#1!(\tikztotarget)!\biglen!-90:(\tikztotarget)$)
    -- ([turn]0,2*\biglen) -- ([turn]0,2*\biglen) -- cycle}},
  half plane/.default={1pt}
}

\newbox{\imageAbox}
\begin{lrbox}{\imageAbox}
    \includegraphics[width=3.5cm, height=3.5cm]{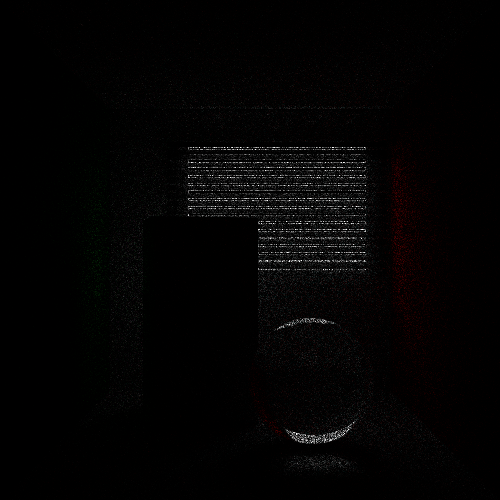} 
\end{lrbox}

\newbox{\imageBbox}
\begin{lrbox}{\imageBbox}
    \includegraphics[width=3.5cm, height=3.5cm]{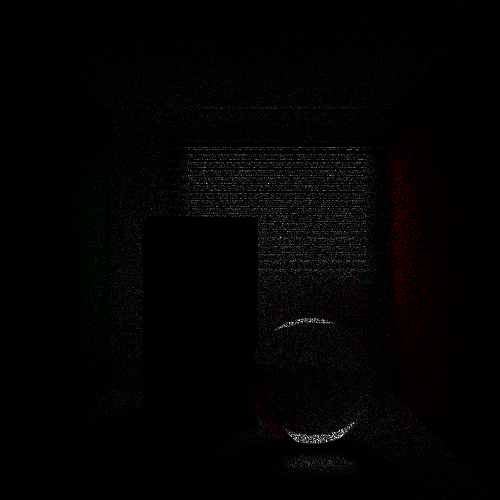}
\end{lrbox}

\newbox{\imageRefAbox}
\begin{lrbox}{\imageRefAbox}
    \includegraphics[width=3.5cm,height=3.5cm]{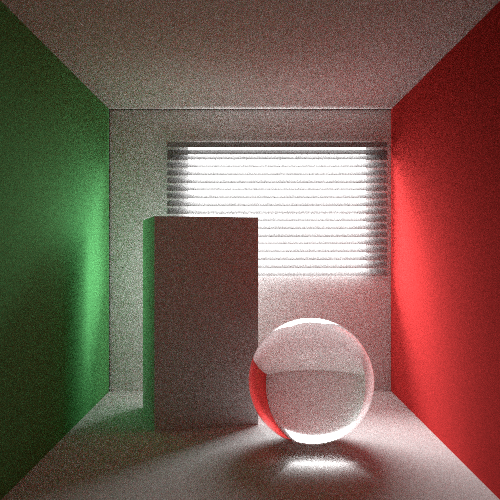} 
\end{lrbox}

\newbox{\imageRefBbox}
\begin{lrbox}{\imageRefBbox}
    \includegraphics[width=3.5cm,height=3.5cm]{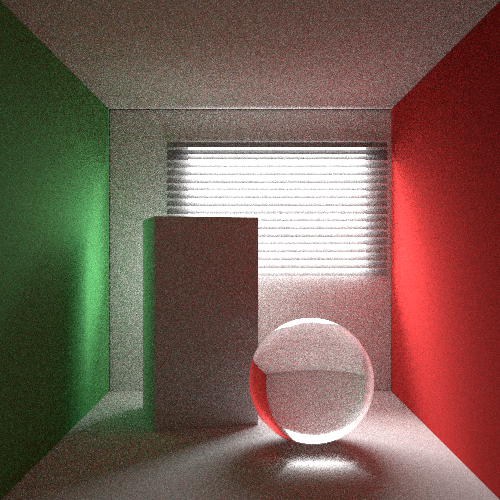} 
\end{lrbox}

\begin{document}

\title{Voronoi integration of the rendering equation}
\date{}

\author[1]{Franck Vandewièle}

\author[2]{Nicolas Chenavier}

\author[1]{Samuel Delepoulle}

\author[1]{Christophe Renaud}

\affil[1]{LISIC, Université du Littoral Côte d'Opale, Calais, France}

\affil[2]{LMPA, Université du Littoral Côte d'Opale, Calais, France}

\affil[ ]{\textit {\{franck.vandewiele, nicolas.chenavier, samuel.delepoulle, christophe.renaud\}@univ-littoral.fr}}


\maketitle

\begin{abstract}In photorealistic image rendering, Monte Carlo methods are the foundation for integration of the rendering equation in modern approaches.
However, despite their effectiveness, traditional Monte Carlo methods often face challenges in controlling variance, resulting in noisy visual artifacts in some difficult to render regions of the image.

In this work, we propose a new approach to the integration of the rendering equation, introducing a Voronoi tessellation reweighting backed by a Poisson point process sampling strategy to address some of the limitations of standard Monte Carlo methods. 

From a theoretical point of view, we show that the variance based on a Poisson-Voronoi tessellation is smaller than that induced by the Monte Carlo method when the intensity of the underlying process is arbitrarily large and the function to be integrated satisfies a Hölder condition.
\end{abstract}


\newpage


\section{Monte Carlo methods}\label{sect:MC_method}

\subsection{Background on Monte Carlo Integration}

Monte Carlo methods are a class of computational algorithms that use random sampling to estimate solutions of integral problems. They do not require prior knowledge of the exact solution and make no demanding assumptions about the functions to integrate. They can be applied to a wide category of problems and are commonly used in various fields such as artificial intelligence, physics, finance, and risk assessment.

Monte Carlo integration has a relatively slow convergence rate, with the error decreasing as $\sfrac{1}{\sqrt{n}}$, where $n$ is the sample size. This means that to halve the error, it is necessary to quadruple the number of samples, which can be computationally expensive for high-precision results \cite{Robert2010}.
 
Nevertheless, it is still a very interesting approach for two main reasons:
\begin{itemize}
 \item the decrease of the error does not depend on the dimension, which makes the method particularly interesting for high-dimensional problems;
 \item the method makes no assumptions about the function other than its integrability.
\end{itemize}
 

Let $f:\RR^d \rightarrow \RR$ a function with compact support integrable over $\RR^d$, equipped with the Euclidean norm $|\cdot |$.

To integrate $f$, Monte Carlo methods assign the same weight to samples $f(x_i)$ taken at locations $(x_i)_{i\leq n}$:


$$ \int_{\RR^d} f(x) \mathrm{d}x \simeq \frac{1}{n} \sum_{i=1}^{n} f(x_i). $$

This approach, by simply averaging the function values, ignores the spatial distribution and density of the samples. In particular, sampling clumps and clusters may be over-represented in the estimation of the integral. 

\subsection{Voronoi reweighting}

\label{sec:voronoireweighting}

Voronoi tessellations \cite{Okabe2000} 
provide an effective representation of the spatial distribution of samples. By defining zones of influence and proximity relations between sites, Voronoi cells offer a convenient and natural means to characterize the samples' geometric features. In particular, reweighting each sample's contribution according to the volume of its Voronoi cell is a straightforward and intuitive method for correcting Monte Carlo estimates, as it naturally mitigates the influence of sample clusters.

Voechovský et al. \cite{Voechovsk2017} apply a Voronoi reweighting scheme to identify under-sampled regions of the integration domain. However, the authors note that simply clipping Voronoi cells at the domain boundaries results in a biased estimation. To ensure an unbiased estimate of the integral, they propose to arrange samples in a periodic tiling, which entails a high computational cost. Guo \& Eisemann \cite{Guo2021}, on the other hand, propose to refine the natural Voronoi reweighting based on clipped cells by applying a correction factor to restore unbiasedness. While this approach is agnostic to the sampling distribution, it is not supported by a theoretical variance upper bound.



In the following, let $S$ be a locally finite subset of 
$\RR^d$,
that is, a set such that its restriction to any compact subset of $\RR^d$ is finite. Denote by
\begin{equation*}
C_{S}(x) = \big\{y\in \RR^d: |y-x|\leq |y'-x| 
\text { for all } y'\in S \setminus\{x\}\big\}
\end{equation*}
the Voronoi cell with nucleus $x$ and by $v_S(x)$ its volume.

We consider this Voronoi weighting:

$$ \sum_{x \in S}f(x) v_S(x), $$ where $f$ is some real function with compact support in $\RR^d$.

Our goal is to analyze the consistency of this estimator and derive a theoretical upper bound on its variance. We demonstrate its effectiveness through experiments on analytical benchmarks as well as in the challenging context of the rendering equation \cite{Kajiya1986}.



\section{Approximation of an integral using a Poisson-Voronoi tessellation}\label{sect:approx}

For theoretical purposes, we consider below a function  $f:\RR^d\rightarrow\RR$ with compact support $W$ satisfying a Hölder condition, namely
\[|f(x)-f(y)|\leq k|x-y|^\alpha,\]
for some $k\geq 0$ and $\alpha \in (0,1]$.
 Let $\eta_n$ be a Poisson point process with intensity $n$ in $\RR^d$, 
i.e. a point process in $\RR^d$ satisfying the two following properties:
\begin{enumerate}[(i)]
    \item for any bounded Borel subset $B$ in $\RR^d$, the random variable $\#(\eta_n \cap B)$ has a Poisson distribution with parameter $n\, \text{vol}_d(B)$;
    \item for any pairwise disjoint Borel subsets $B_1,\ldots, B_k$, $k\geq 1$, in $\RR^d$, the random variables $\#(\eta_n\cap B_1)$, \ldots $\#(\eta_n\cap B_k)$ are independent.
\end{enumerate}
From a theoretical point of view, such a process is natural because it allows expectations of sums to be computed using Slivnyak-Mecke  formula (see Appendix \ref{annex:proof}). In practice, to simulate it in a bounded window $W$, one simply generates a random number 
$N$ following a Poisson distribution with parameter $n\, \text{vol}_d(W)$ and then, conditionally on $N$, throws 
$N$ points independently of each other and uniformly at random in 
$W$. Using notations similar to those in Section \ref{sec:voronoireweighting}, we denote by
\begin{equation*}
C_{\eta_n}(x) = \big\{y\in \RR^d: |y-x|\leq |y'-x| 
\text { for all } y'\in \eta_n\setminus\{x\}\big\}        
\end{equation*}
 the Voronoi cell with nucleus $x\in \eta_n$ and by $v_{\eta_n}(x)$ its volume.  

\begin{Prop}
\label{prop:variance}
With the above notations, we have
\begin{enumerate}[(i)]
\item $\EEE{\sum_{x\in \eta_n}f(x)v_{\eta_n}(x)} = \int_{\RR^d}f(x)\mathrm{d}x$;
\item $\VVV{\sum_{x\in \eta_n}f(x)v_{\eta_n}(x)} \leq c n^{-1-\tfrac{2\alpha}{d}}$, 
for some $c>0$ depending only on $d$ and $f$.
\end{enumerate}
\end{Prop}

Approximating the integral of a given function (especially when it is Hölder continuous) using the values it takes on a set of points, weighted by the volume of Voronoi cells, is classical and naturally arises in quantization problems (see \cite{Pages98}). However, the underlying methods consider a deterministic set of points rather than a (random) point process and are computationally expensive.

Using a Voronoi tessellation based on a Poisson point process on $\RR^d$ (and not only on the support of $f$) is important in order to obtain assertion \textit{(i)}, because the underlying process is stationary. Conversely, a Voronoi tessellation based on a point process defined only on the support of $f$ would yield a biased estimator of the integral, since it would not be stationary.

Concerning (ii), we notice that the upper bound, namely $ n^{-1-\tfrac{2\alpha}{d}}$, is the term of a convergent series. In particular, thanks to the the Tchebychev's inequality and the Borel-Cantelli lemma, it implies 
\[\sum_{x\in \eta_n}f(x)v_{\eta_n}(x) \conv[n]{\infty}  \int_{\RR^d}f(x)\mathrm{d}x \quad \text{ a.s..}\]

It is noteworthy that the variance converges to zero at a faster rate than that of the Monte Carlo method, for which the rate is $O(n^{-1})$. To the best of our knowledge, such a result is new. 

An interesting point is that, when $f$ does not satisfy a Hölder condition but is of the form $f=\ind{W}$, where $W$ is some admissible set of $\RR^d$, then
\[\VVV{\sum_{x\in \eta_n}f(x)v_{\eta_n}(x)}  \leq c n^{-1-\tfrac{2}{d}} = O(n^{-1-\tfrac{1}{d}});\]
see \cite{Heveling09, Schulte12bis,Thale16}. See also \cite{Reitzner24} for a recent work on Poisson-Delaunay approximation. Here again, our rate of convergence decreases faster.

See appendix \ref{annex:proof} for a proof of this result.


\section{Application to integration problems}

\subsection{Leveraging Poisson-Voronoi tessellations for numerical estimations}
\label{section:margin}
Proposition \ref{prop:variance} offers an appealing reduction technique for Monte Carlo estimations of Hölder continuous functions. In the following, we will consider the numerical computation of the integral of real functions defined in $\RR^2$ with compact support $\WP=\WPP$. 
To apply this result to this kind of problem, we need to sample a real function with a Poisson point process in $\RR^2$. To restrict a Poisson process $\eta_n$ to a finite sampling of a compact window $\WP$, we first draw an integer $N$ from a Poisson distribution with parameter $|\WP|n = n$. Then, given $N$, $x_1$, \ldots, $x_N$ are sampled uniformly from $\WP$.

Even if this construction allows us to reduce it to a sampling equivalent to a Poisson process, our estimation still faces a difficulty: a finite set of nuclei $x_1, \ldots, x_N$ gives rise to unbounded Voronoi cells.

For this reason, we consider a \emph{Stretched window Poisson Point Process} (SPPP) sampling strategy: instead of considering a Voronoi tessellation from nuclei drawn in $\WP=\WPP$, we draw random nuclei from a stretched window $W' = \WPPe$, where  $\eps > 0$ is chosen to ensure that Voronoi cells of the nuclei contained in $\WP$ are bounded, with a certain level of confidence.

Note that the extra nuclei drawn from outside $\WP$ are used to constrain the Voronoi tessellation, and not for sampling function $f$. By choosing a suitable confidence level, the acceptance-rejection method yield a valid Voronoi partition for estimating the integral of $f$ over $\WP$ via the result of Proposition \ref{prop:variance} with minimum overhead.

To choose $\eps$, we rely on the following argument. 

Let us assume that $n$ points are distributed within a window $\WP = \WPP$. Let $\WPe = \WPPe$ the corresponding stretched window for some $\eps$. We process \emph{as if} dealing with a Poisson point process of intensity $n$ in the plane. Let $X_n = \{x_1, \ldots, x_n\}$  denote $n$ points sampled independently and uniformly at random in
$\WP$, and $X'_n = \{x'_1, \ldots x'_{n'}\}$ denote the $n'$ points (also thrown independently and uniformly at random) in the strip $\WPe \setminus \WP$, distributed independently of the points in $X$. Since the area of the strip $\WP' \setminus \WP$ is $4\eps + 4\eps^2$, if we want to preserve the density, the number of points $n'$ distributed in the strip is given by $n' = \lfloor n(4\eps + 4\eps^2)\rfloor$. Let ${\Tilde{X}_n = X_n \cup X'_n}$ and observe that ${\Tilde{X}_n \cap \WP = X_n}$. Recall that $\eta_n$ denote a Poisson point process of intensity $n$. We have :
\begin{align*}
    \PPP{\text{reject the simulation} }   & \quad = \PPP{\exists x \in \Tilde{X}_n\cap \WP: C_{\Tilde{X}_n}(x) \text{ is unbounded} } \\
    & \quad \leq \PPP{\exists x \in \Tilde{X}_n\cap \WP: R(C_{\Tilde{X}_n}(x)) > \eps} \\
    & \quad \eq[n]{\infty} \PPP{\exists x \in \eta_n\cap \WP: R(C_{\eta_n}(x)) > \eps}
\end{align*}
where, for any locally finite subset $\chi\subset \RR^2$, the quantity $R_\chi(x)$ denotes the circumscribed radius associated with $C_\chi(x)$, i.e.
\[R_\chi(x) = \inf\{r>0: B(x,r) \supset C_\chi(x)\}.\] 

Moreover
\begin{align*}
  \PPP{\exists x \in \eta_n\cap \WP: R(C_{\eta_n}(x)) > \eps}  &\quad\leq \EEE{\sum_{x\in \eta_n\cap \WP}\ind{R(C_{\eta_n}(x))>\eps}}\\
  &\quad =  n\int_{\WP}\PPP{R(C_{\eta_n\cup\{x\}}(x))>\eps}\mathrm{d}x\\
  &\quad = n\PPP{R(C_{\eta_n\cup\{0\}}(0))>\eps}\\
  &\quad = n\PPP{R(C_{\eta_1\cup\{0\}}(0))>\eps n^{1/2}}.
\end{align*}
In the above expression, the first equality follows from the Slivnyak-Mecke formula (see e.g. Corollary 3.2.3 in \cite{Schneider08}), the second equality follows from the facts that the area of $\WP$ equals 1 and that $\eta_n$ is stationary, and  the third one follows from the relation $\eta_n \overset{\mathcal{D}}{=} n^{-1/2}\eta_1$. Thanks to \cite{Calka02bis}, Eq. (6) this gives
\[\PPP{\exists x \in \eta_n\cap \WP: R(C_{\eta_n}(x)) > \eps} \leq 4\pi 
n^2\eps^2 e^{-\pi n\eps^2}.\] 
Therefore, as $n$ is large, the probability to reject the simulation is lower than $4\pi n^2\eps^2 e^{-\pi n\eps^2}$.


In our experimental setup, we use a target level of confidence of $\frac{1}{1,000}$. With $n=10,000$, for instance, this gives $\eps \simeq 0.0255$.

\subsection{Poisson-Voronoi estimators}


We recall the notation introduced previously. $X_n = \{x_1, \ldots, x_n\}$ represents a set of $n$ points sampled independently and uniformly in the compact window $\WP = \WPP$; $\WPe$ denotes the stretched window $\WPPe$ for some $\eps > 0$. Conversely, the set of supplementary points ${X'_n} = \{x'_1, \ldots x'_{n'}\}$ is sampled independently and uniformly in the strip $\WP' \setminus \WP$, and $\Tilde{X}_n = X_n \cup X'_n$ is the augmented set of initial points.

\begin{figure}[htbp]
    \centering
    \subfloat[Voronoi estimator $E_V$ takes every cell with nucleus inside $\WP$ into account.\label{fig:pixel1}]{\includegraphics[width=.33\textwidth]{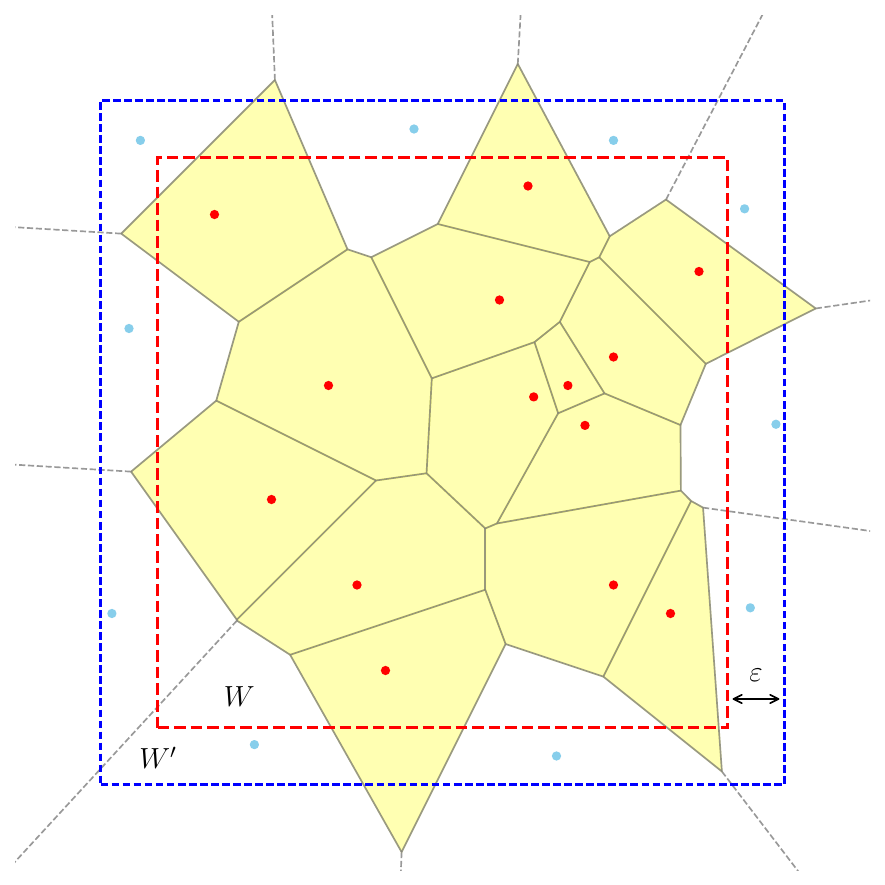}}\hspace*{2em}
    \subfloat[Filtered Voronoi estimator $E_F$ discards cells extending beyond $\WPe$.\label{fig:pixel2}]{\includegraphics[width=.33\textwidth]{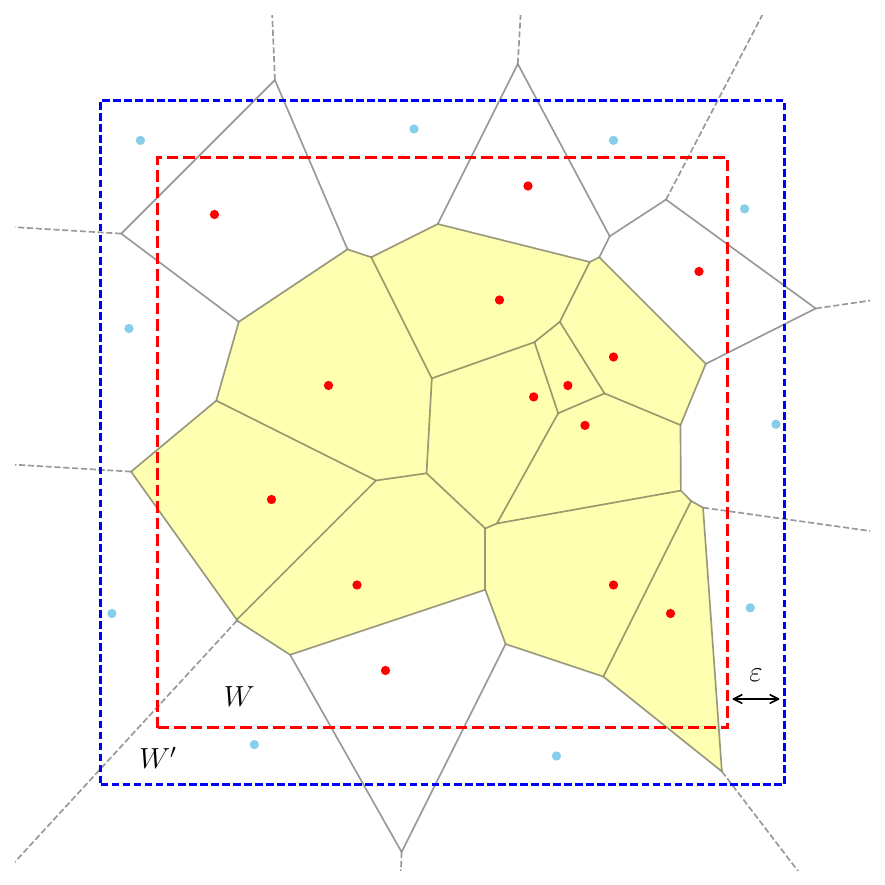}}
    
    \subfloat[Clipped Voronoi estimator $E_C$ clips cells to the boundaries of $\WP$.\label{fig:pixel3}]{\includegraphics[width=.33\textwidth]{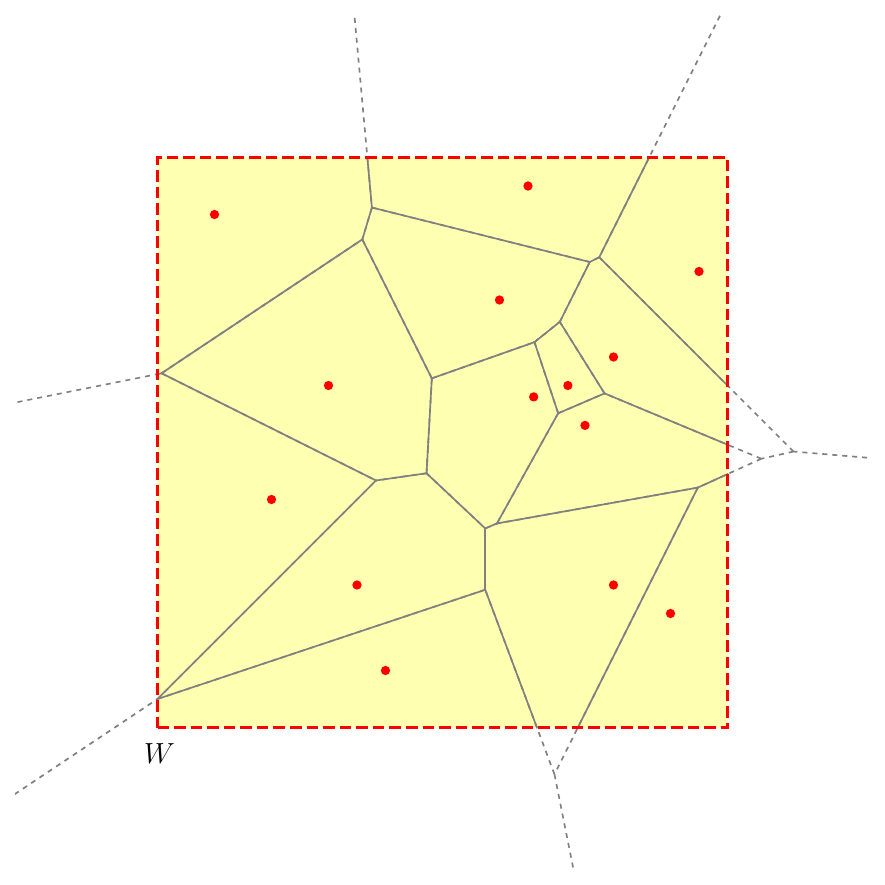}}%
    \caption{Three Voronoi integration estimators over $\WP=\WPP$. Voronoi cells nuclei inside $\WP$ are depicted in red, and auxiliary nuclei in blue. Cells used by the estimators are highlighted in yellow.
    }
    \label{fig:pixel_distribution}
\end{figure}

Let $f : \RR^2 \rightarrow \RR$ be a function with compact support $\WP$. From proposition \ref{prop:variance}, we derive three estimators for integration problems. Figure \ref{fig:pixel_distribution} summarizes the main characteristics of each estimator. 

\begin{description}
\item[Voronoi estimator] is the immediate adaptation of proposition \ref{prop:variance} for a SPPP (Figure \ref{fig:pixel1}):
$$E_V (f) = \sum_{x\in {{X}_n}}f(x)v_{\Tilde{X}_n}(x).$$

\item[Filtered Voronoi estimator] is another take on proposition \ref{prop:variance}. While SPPP guarantees that the Voronoi cells are bounded, in practice, their diameter might be very large. For this reason, this estimator filters out Voronoi cells that might extend beyond $\WPe$. Let $X^{-}_n \subseteq X_n$ denote the subset of nuclei for which the associated Voronoi cell, computed based on the augmented configuration $\Tilde{X}_n$, lies entirely within the domain $W'$ (Figure \ref{fig:pixel2}). We define the filtered Voronoi estimator the following way: 
\newcommand{\XXX}{X^{-}_n}
$$ E_F(f) = \frac{1}{\sum_{x\in \XXX}  v_{\Tilde{X}_n}(x)} \sum_{x\in \XXX} f(x) v_{\Tilde{X}_{n}}(x). $$
\item[Clipped Voronoi estimator] discards the stretching strategy entirely and clips Voronoi cells on the borders of $\WP$ (Figure \ref{fig:pixel3}). This estimator is expected to be biased:

$$E_C(f) = \sum_{x\in {X_n}}f(x)|C_{X_{n}}(x)\cap \WP|. $$
\end{description}



\section{Experiment}

\subsection{Experimental goals and hypothesis}

The primary goal is to observe that the Voronoi weighting approach, which uses the area of Voronoi cells to weight sample contributions, provides a more efficient and stable estimate of an integral compared to the unweighted average used in Monte Carlo. This efficiency and stability are quantified by comparing the standard deviation of the estimates for different integration methods.

A list of functions is used to test the convergence of the methods (see table \ref{table:function_list} for their definition).

Each of these functions exhibit some form of irregularity at the origin:

\begin{itemize}
    \item \texttt{holder\_\textalpha} is Hölder continuous with exponent $\alpha\in (0, 1]$; it is not differentiable at the origin if $\alpha < 1$;
    \item \texttt{not\_holder} is continuous on $W$ but not Hölder  and therefore not differentiable (at the origin); 
    \item \texttt{discontinuity} is undefined at the origin, thus not Hölder continuous on $\WP$, yet $\mathcal{C}^\infty$ on $\WP \setminus \{(0, 0)\}$.
\end{itemize}


\begin{table}[htbp]
    \centering
    \caption{Definition of the used functions}
    \label{table:function_list}
    \begin{tabular}{l p{0.6\columnwidth}}
        \toprule
        \textbf{name} & \textbf{definition} \\
        \midrule
        \texttt{holder\_\textalpha} & $| x  y |^\alpha$ \\
        \addlinespace
        \texttt{not\_holder} & $x y \sin\left(\frac{1}{x}\right)\sin\left(\frac{1}{y}\right)$ \\
        \addlinespace
         \texttt{discontinuity} & $\left( x^2 + y^2\right)^{-\frac{1}{2}}$ \\
        \bottomrule
    \end{tabular}
\end{table}


\begin{table*}[ht]
    \centering
    \caption{Comparison of standard deviations and average execution durations (ms) over $10,000$ runs at $4,096$ samples per pixel (spp).}
    \label{tab:comparison_no_nvor}
\begin{tabular}{|l||c|c|c||c|c|c|}
    \hline
    \multicolumn{1}{|l||}{} & \multicolumn{3}{c||}{\textbf{Standard Deviation}} & \multicolumn{3}{c|}{\textbf{Time}} \\
    \hline
    \textbf{Fonction} & \textbf{\texttt{mc}} & \textbf{\texttt{vor}} & \textbf{\texttt{fvor}} & \textbf{\texttt{mc}} & \textbf{\texttt{vor}} & \textbf{\texttt{fvor}} \\
    \hline
    \hline
    \texttt{holder\_1} & 0.00086 & 0.000171 & \textbf{0.000167} & \textbf{2.280} & 201.970 & 149.870 \\
    \texttt{holder\_0.5} & 0.001786 & \textbf{0.000291} & 0.000292 & \textbf{4.040} & 204.820 & 150.300 \\
    \texttt{holder\_0.1} & 0.001461 & 0.000299 & \textbf{0.000297} & \textbf{3.780} & 216.570 & 156.930 \\
    \texttt{holder\_0.01} & 0.000244 & $\mathbf{<10^{-6}}$ & $\mathbf{<10^{-6}}$ & \textbf{4.070} & 229.090 & 151.750 \\
    \texttt{not\_holder} & 0.000585 & 0.000179 & \textbf{0.000178} & \textbf{5.450} & 228.610 & 159.330 \\
    \texttt{discontinuity} & 0.10716 & \textbf{0.079052} & 0.10411 & \textbf{2.530} & 231.970 & 152.400 \\
    \hline
\end{tabular}
\label{tab:comparison}
\end{table*}


\subsection{Results}

\subsubsection{Experimental labeling} In the following experiments, the standard Monte Carlo estimation algorithm is denoted \texttt{mc}. Voronoi estimators $E_V$, $E_F$ and $E_C$ are denoted \texttt{vor}, \texttt{fvor} and \texttt{cvor} respectively.




\begin{figure}[htbp]
     \centering
     \includegraphics[width=0.65\columnwidth]{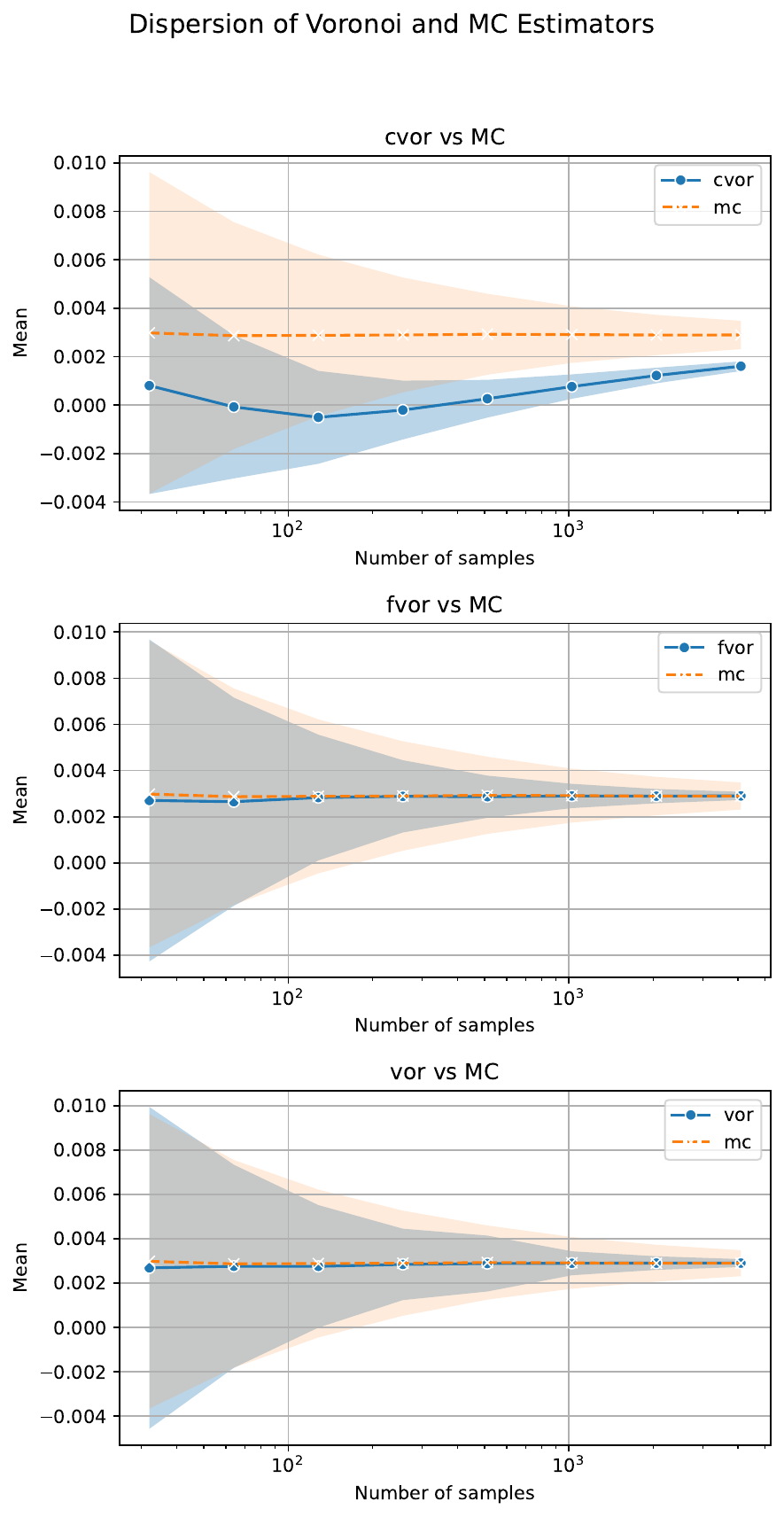}
     \caption{Performance of integrators for function \texttt{not\_holder}. Each plot compares a reweighting integration method against standard Monte Carlo (labeled \texttt{mc}). Solid lines represent the average estimate over 10,000 runs, with sample count $n$ ranging from $2^5$ to $2^{14}$. Shaded regions indicate the dispersion of these estimates within $\pm 1$ standard deviation. The estimators $E_V$, $E_C$ and $E_F$ (labeled \texttt{vor}, \texttt{cvor} and \texttt{fvor} respectively) demonstrate significant variance reduction. As expected, the clipped Voronoi reweighting method \texttt{cvor} exhibits estimation bias, while still reducing the variance.}
     \label{fig:variance_convergence}
\end{figure}

\begin{figure}[htbp]
     \centering
     \includegraphics[width=0.95\columnwidth]{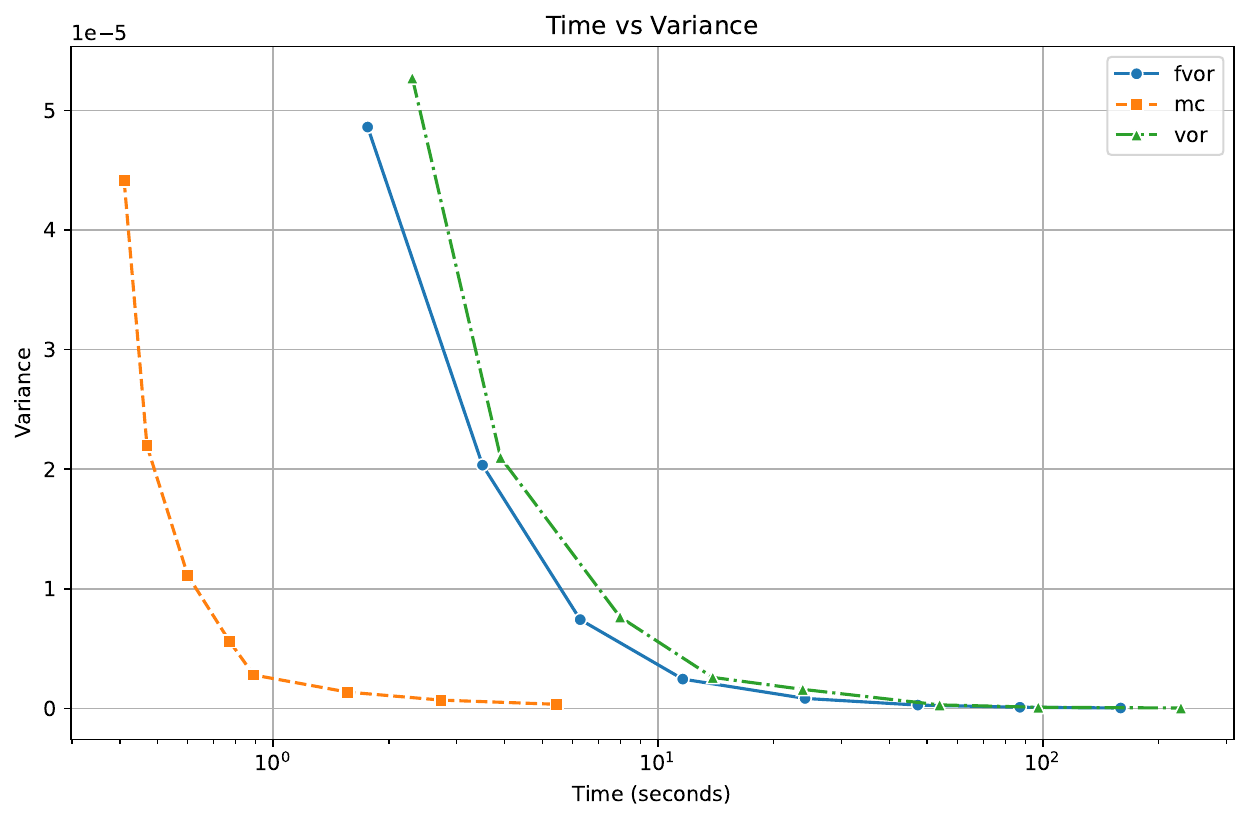}
     \caption{Comparative efficiency of three algorithms for evaluating the integral of the function \texttt{not\_holder}. 
     }
     \label{fig:link_bench_not_c1}
\end{figure}     


\subsubsection{Analysis}
Figure \ref{fig:variance_convergence} shows that all three estimators achieve significant variance reduction for function \texttt{not\_c1}. The clipped estimator $E_C$, however, exhibits significant estimation bias. Consequently, being of limited practical interest, this estimator will not be systematically evaluated in subsequent experiments.
Table \ref{tab:comparison} shows that $E_V$ and $E_F$ achieve significant variance reduction for the Hölder function cases considered. Though Proposition \ref{prop:variance} does not hold for this case, a comparable variance reduction is also observed for the function \texttt{not\_holder}. 
However, Voronoi estimators fail to outperform classical Monte Carlo for function \texttt{discontinuity}, which has a singularity at the origin. 




During our experiments, we found that in practice, the SPPP sampling strategy leads to a much lower rejection rate than implied by the target confidence level.

In addition to generating the random sample distribution, standard Monte Carlo requires only repeated evaluations of the integrand. For explicitly defined functions, this is a computationally trivial task. In contrast, our method requires computing a Voronoi tesselation, a costly operation that dominates the total runtime. Table \ref{tab:comparison} and Figure \ref{fig:link_bench_not_c1} show that execution times disqualify Voronoi reweighting methods for the integration of explicitly defined functions. 
However, for functions where evaluation is costly, and specifically of an order of magnitude comparable to the Voronoi tessellation construction, the rapid convergence of our method might mitigate the computational overhead. For this reason, the next section is dedicated to applying this approach to the rendering equation.

\section{Results with MC integration for Global illumination}




\subsection{The Rendering Equation and its Monte Carlo Estimation}
\begin{figure}[htbp]
		\tdplotsetmaincoords{70}{110}
	\centering\begin{tikzpicture}[
		tdplot_main_coords,
		scale=1.5,
		>=Stealth,
		pixel/.style={thin, draw=gray!50},
		highlight/.style={thick, dashed, draw=red!95!black, fill=blue!10, fill opacity=0.3}
		]
		
		\coordinate (C) at (0,0,0); 
		\def\focal{4} 
		
		\def\py{0} 
		\def\pz{0} 
		\def\sz{1}   
		
		\coordinate (P1) at (\focal, \py, \pz);         
		\coordinate (P2) at (\focal, \py+\sz, \pz);     
		\coordinate (P3) at (\focal, \py+\sz, \pz+\sz); 
		\coordinate (P4) at (\focal, \py, \pz+\sz);     
		
		\coordinate (S) at (\focal, \py+0.45*\sz, \pz+0.7*\sz);
		\coordinate (S1) at (\focal, \py+0.2*\sz, \pz+0.1*\sz);
		\coordinate (S2) at (\focal, \py+0.1*\sz, \pz+0.9*\sz);
		\coordinate (S3) at (\focal, \py+0.9*\sz, \pz+0.3*\sz);
		\coordinate (S4) at (\focal, \py+0.65*\sz, \pz+0.2*\sz);

		\shade[left color=orange!20, right color=orange!5, opacity=0.4] 
		(C) -- (P1) -- (P4) -- cycle; 
		\shade[left color=orange!20, right color=orange!5, opacity=0.4] 
		(C) -- (P4) -- (P3) -- cycle; 
		\shade[left color=orange!30, right color=orange!10, opacity=0.4] 
		(C) -- (P3) -- (P2) -- cycle; 
		\shade[left color=orange!30, right color=orange!10, opacity=0.4] 
		(C) -- (P2) -- (P1) -- cycle; 
		
		\draw[orange!80!black, thin, dashed] (C) -- (P1);
		\draw[orange!80!black, thin] (C) -- (P2);
		\draw[orange!80!black, thin] (C) -- (P3);
		\draw[orange!80!black, thin] (C) -- (P4);

		\begin{scope}[canvas is yz plane at x=\focal]
			\draw[pixel, step=1] (-2, -2) grid (2, 2);
			
			\draw[thick, gray] (-2, -2) rectangle (2, 2);
			
			\node[anchor=south west, black] at (-2, -2) {Pixel grid};
		\end{scope}

		\draw[highlight] (P1) -- (P2) -- (P3) -- (P4) -- cycle;

		\coordinate (S_far) at ($ (C)!1.6!(S) $);
		\coordinate (S1_far) at ($ (C)!1.6!(S1) $);
		\coordinate (S2_far) at ($ (C)!1.6!(S2) $);
		\coordinate (S3_far) at ($ (C)!1.6!(S3) $);
		\coordinate (S4_far) at ($ (C)!1.6!(S4) $);
		\draw[orange!80!black, dashed, -] (C) -- (S) node[above, pos=1.0] {};
		\draw[orange!80!black, thick, ->] (S) -- (S_far) node[above, pos=1.0] {};
		\draw[orange!80!black, dashed, -] (C) -- (S1) node[above, pos=1.0] {};
		\draw[orange!80!black, thick, ->] (S1) -- (S1_far) node[above, pos=1.0] {};
		\draw[orange!80!black, dashed, -] (C) -- (S2) node[above, pos=1.0] {};
		\draw[orange!80!black, thick, ->] (S2) -- (S2_far) node[above, pos=1.0] {};
		\draw[orange!80!black, dashed, -] (C) -- (S3) node[above, pos=1.0] {};
		\draw[orange!80!black, thick, ->] (S3) -- (S3_far) node[above, pos=1.0] {};
		\draw[orange!80!black, dashed, -] (C) -- (S4) node[above, pos=1.0] {};
		\draw[orange!80!black, thick, ->] (S4) -- (S4_far) node[above, pos=1.0] {};
		\draw[black] (S2_far) node[left=3pt] {$L_1$};
		\draw[black] (S1_far) node[left=1pt, below=1pt] {$L_2$};
		\draw[black] (S_far) node[left=3pt] {$L_3$};
		\draw[black] (S3_far) node[below=1pt] {$L_4$};
		\draw[black] (S4_far) node[below=1pt] {$L_5$};
		\fill[red] (S) circle (1.8pt);
		\fill[red!95!black] (S1) circle (1.8pt);
		\fill[red!90!black] (S2) circle (1.8pt);
		\fill[red!80!black] (S3) circle (1.8pt);
		\fill[red!70!black] (S4) circle (1.8pt);
		
		\node[anchor=south west, black]  at ([yshift=-2mm] P2) {$W$};
		
		\fill[black] (C) circle (1.5pt) node[right=5pt] {$C$};
	\end{tikzpicture}
    \caption{\label{fig:mc-contributions-pixel}Simulating numerous potential light paths connecting a virtual camera $C$ to light sources leads to the accumulation of light contributions $(L_i)$ over a pixel area $W$, defining its final color. In a standard Monte Carlo framework, these contributions are averaged without accounting for the geometry of their distribution.}
\end{figure}
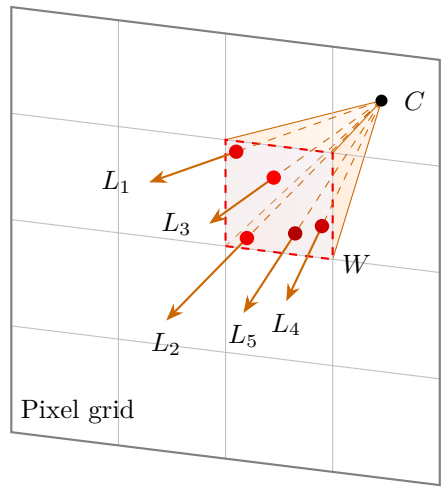

Path tracing is a Monte Carlo estimator for the rendering equation: pixel values are expressed by an integral which is numerically solved by averaging random samples of light transport paths.

The rendering equation \cite{Kajiya1986} is the fundamental mathematical model that describes how light interacts with surfaces in a scene. For any non trivial scene, this equation has no analytic solution. To estimate its solution, the Path Tracing algorithm relies on a stochastic simulation of the interaction of light with the geometry and the materials of the scene. Being a stochastic algorithm, Path Tracing depends on variance reduction techniques for achieving an accurate, noise-free visual result. Figure \ref{fig:mc-contributions-pixel} illustrates how the color of each individual pixel to display on the screen is computed: by simulating numerous light paths traversing it \cite{Pharr2023}.

\subsection{Experimental renderer}

YAPT stands for Yet Another Path Tracer \cite{Vandewiele2025}. It is an experimental rendering engine developed to investigate and test novel sampling and integration techniques in the field of physically-based rendering. As a path tracer, YAPT simulates the way light interacts with objects in a virtual environment to produce  realistic images (Figure \ref{fig:trois_images}).

The core architectural principle is the separation between the two fundamental processes of Monte Carlo integration:

\begin{description}
\item[Sampler] responsible for generating the sample distribution used to estimate the light transport integral. This component determines the spatial distribution of light samples (e.g., uniform, stratified, or SPPP)

\item[Integrator] responsible for calculating the light transport (radiance) along the paths defined by the sampler. This component implements the specific rendering algorithm, such as classic Path Tracing or techniques incorporating explicit light sampling like Next Event Estimation (NEE), which is used in these experiments.
\end{description}

\subsection{Convergence and Efficiency Analysis}

\begin{figure}[htbp]
\centering
\subfloat[16 samples per pixel\label{fig:image1}]{\includegraphics[width=.28\textwidth]{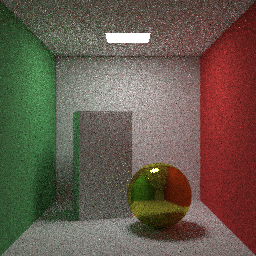}}
\hfill
\subfloat[256 samples per pixel\label{fig:image2}]{\includegraphics[width=.28\textwidth]{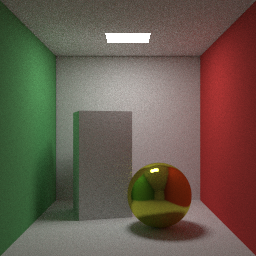}}
\hfill
\subfloat[$4,096$ samples per pixel\label{fig:image3}]{\includegraphics[width=.28\textwidth]{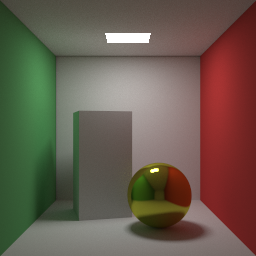}}%
\caption{Illustration of Monte Carlo rendering convergence with YAPT. The estimate gradually approaches the exact scene radiance as the sample count increases, resulting in a progressively more accurate estimation of the scene's true radiance.}
\label{fig:trois_images}
\end{figure}

The core purpose of YAPT is evaluating the performance and convergence of new integration methods compared to standard Monte Carlo techniques. This engine was used to evaluate the image rendering performance of estimators $E_V$ and $E_F$, for which experimental results in section 4 demonstrate a variance reduction for Hölder functions integrated over a compact window. In the context of the rendering equation, such an assumption is in fact not very restrictive, since pixels are squares, hence bounded sets, and the radiance can be modeled by a sufficiently regular function.
Experiments were performed using the classic Cornell Box scene, a standard benchmark in computer graphics. Figure \ref{fig:trois_images} visually illustrates the convergence process inherent to Monte Carlo rendering. By progressively increasing the sample count per pixel, the approximated image approaches the true solution, resulting in a reduction of variance, perceived as visual noise. During our experiments, the following rendering settings were chosen : sample size ranging from 64 to 4,096, maximum path depth 8, and Next Event Estimation. 



\subsection{Quantifying Estimator Performance}

\begin{figure}[htbp]
    \centering
    \subfloat[Voronoi estimation versus uniform Monte Carlo\label{fig:graph1-spp-vor}]{\includegraphics[width=.47\textwidth]{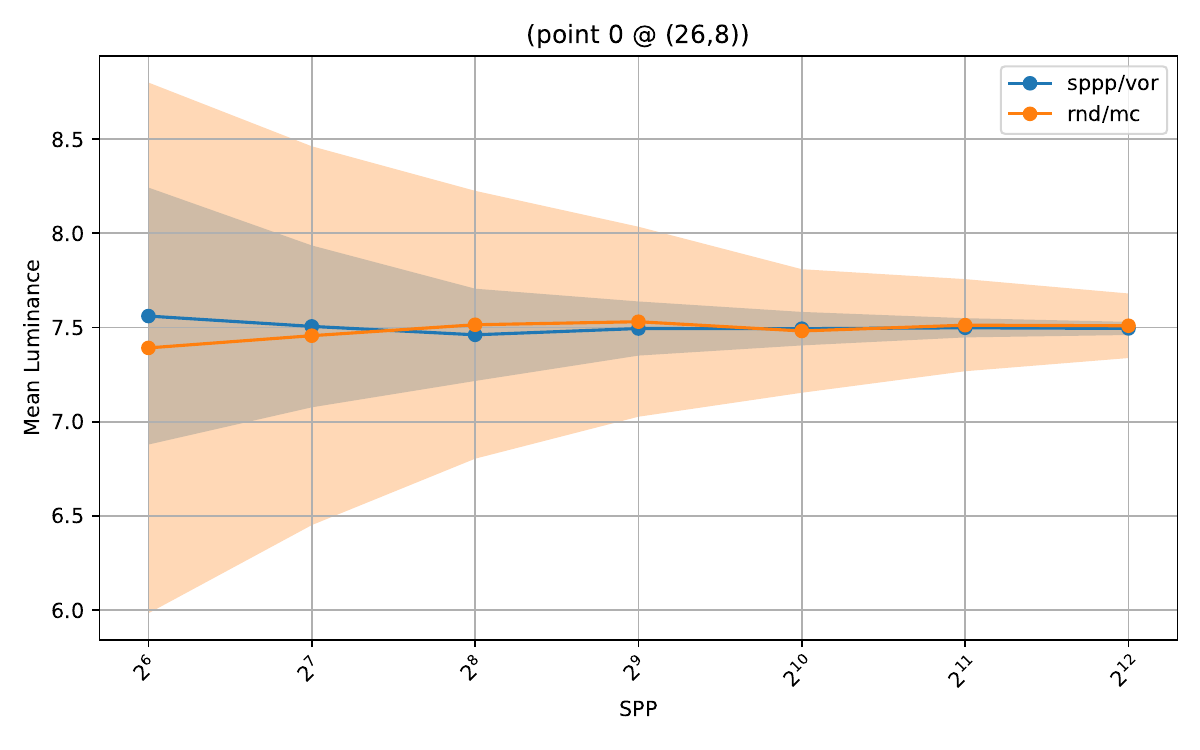}}\hfill
        \subfloat[Filtered Voronoi estimation vs Monte Carlo\label{fig:graph2-spp-fvor}]{\includegraphics[width=.47\textwidth]{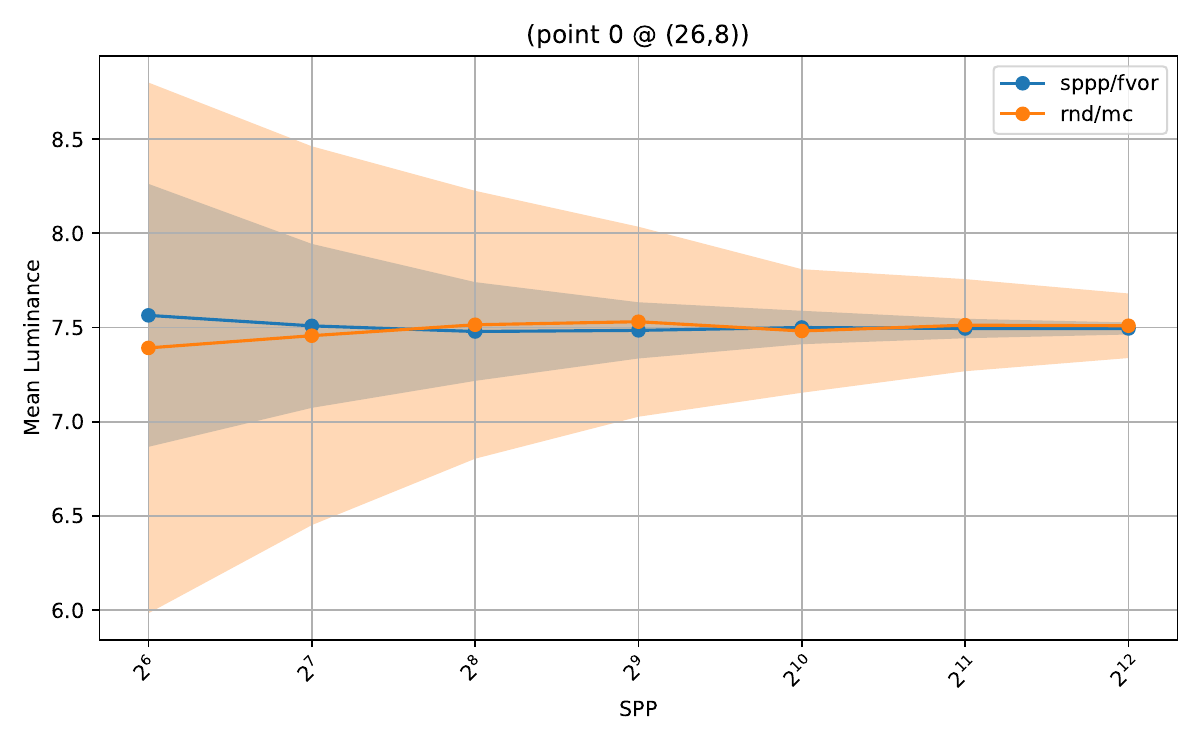}}\hfill
    \caption{Evolution of the mean radiance estimate and $\pm 1$ standard deviation interval as a function of the sample size. Results obtained with the Voronoi estimator $E_V$ \texttt{sppp/vor} and the filtered Voronoi estimator $E_F$ \texttt{sppp/fvor} are compared to those from random sampling using Monte Carlo estimation  \texttt{rnd/mc} in each graph. For this pixel, located along a source boundary, both Voronoi estimators exhibit a faster convergence rate than standard MC estimation.}
    \label{fig:yapt-convergence}
\end{figure}

Figure \ref{fig:yapt-convergence} provides a quantitative analysis of different estimators by plotting the mean estimated radiance. Shaded regions represent the $\pm 1$ standard deviation interval.
These evaluations were performed for a specific, difficult-to-sample pixel located near a light source boundary.
Accurately evaluating the color of such a pixel is indeed difficult. During the sampling process, contributions collected from direct lighting in the source surface will have high amplitude, whereas off-source contributions, which are indirectly lit, will be evaluated with significant variance, due to the diffuse nature of the material of the walls of the scene.



Figure \ref{fig:yapt-convergence} demonstrates that Voronoi integration methods achieve lower variance than standard Monte Carlo approach for identical sample sizes when target pixel is located at the boundary of a light source.

\begin{figure*}[htbp]
\centering
    \newlength{\gapsize}
    \setlength{\gapsize}{6mm}
    \setlength{\tabcolsep}{0pt}
    \renewcommand{\arraystretch}{0}
    \begin{tabular}{c@{\hspace{\gapsize}}c}
        \subfloat[regular Monte Carlo with sample size 256\label{fig:refA}]{%
            \begin{tikzpicture}[baseline=(current bounding box.north)]
                \node[inner sep=0pt] at (0,0) {\usebox{\imageRefAbox}};
            \end{tikzpicture}%
        } &
        
        \subfloat[difference map with reference image (MSE = 0.030)\label{fig:imgA}]{%
            \begin{tikzpicture}[
                spy using outlines={rectangle, red, magnification=5, size=3.5cm, connect spies}, 
                baseline=(current bounding box.north)
            ]
                \node[inner sep=0pt] (main) at (0,0) {\usebox{\imageAbox}};
                
                \spy on (0.0, 0.5) in node at (3.5cm + \gapsize, 0); 
            \end{tikzpicture}%
        } \\
        
        \vspace{1mm} & \vspace{1mm} \\

        
        \subfloat[Voronoi ponderation \texttt{vor} with sample size 256\label{fig:refB}]{%
            \begin{tikzpicture}[baseline=(current bounding box.north)]
                \node[inner sep=0pt] at (0,0) {\usebox{\imageRefBbox}};
            \end{tikzpicture}%
        } &
        
        \subfloat[difference map with reference image (MSE = 0.015)\label{fig:imgB}]{%
            \begin{tikzpicture}[
                spy using outlines={rectangle, blue, magnification=5, size=3.5cm, connect spies}, 
                baseline=(current bounding box.north)
            ]
                \node[inner sep=0pt] (main) at (0,0) {\usebox{\imageBbox}};
                
                \spy on (0.0, 0.5) in node at (3.5cm + \gapsize, 0); 
            \end{tikzpicture}%
        } \\
        
    \end{tabular}
    \caption{This figure compares two integration methods in a Path Tracer for a fixed number of samples per pixel of $n=256$.
    Top row: (a) Rendering of the Cornell Box scene using random Monte Carlo (MC) integration for $n=256$ samples per pixel (spp). (b) Error map (difference) between the MC rendering (a) and a high-quality reference (reference calculated with $n =4096$ spp and stratified sampling).
    Bottom row: (c) Rendering of the same scene using Voronoi-weighted integration for $n=256$ spp. (d) Map of absolute errors (difference) between the Voronoi rendering (c) and the same reference image.
    A comparison of error maps (b) and (d) demonstrates that Voronoi-weighted integration \texttt{vor} (c) yields lower variance than standard MC integration at equal sample counts, most notably at radiance discontinuities, such as shadows and light source edges.}
    \label{fig:main_figure}
\end{figure*}



\begin{figure}[htbp]
    \centering
    \subfloat[Near a light source\label{fig:detail1}]{\includegraphics[width=.4\textwidth]{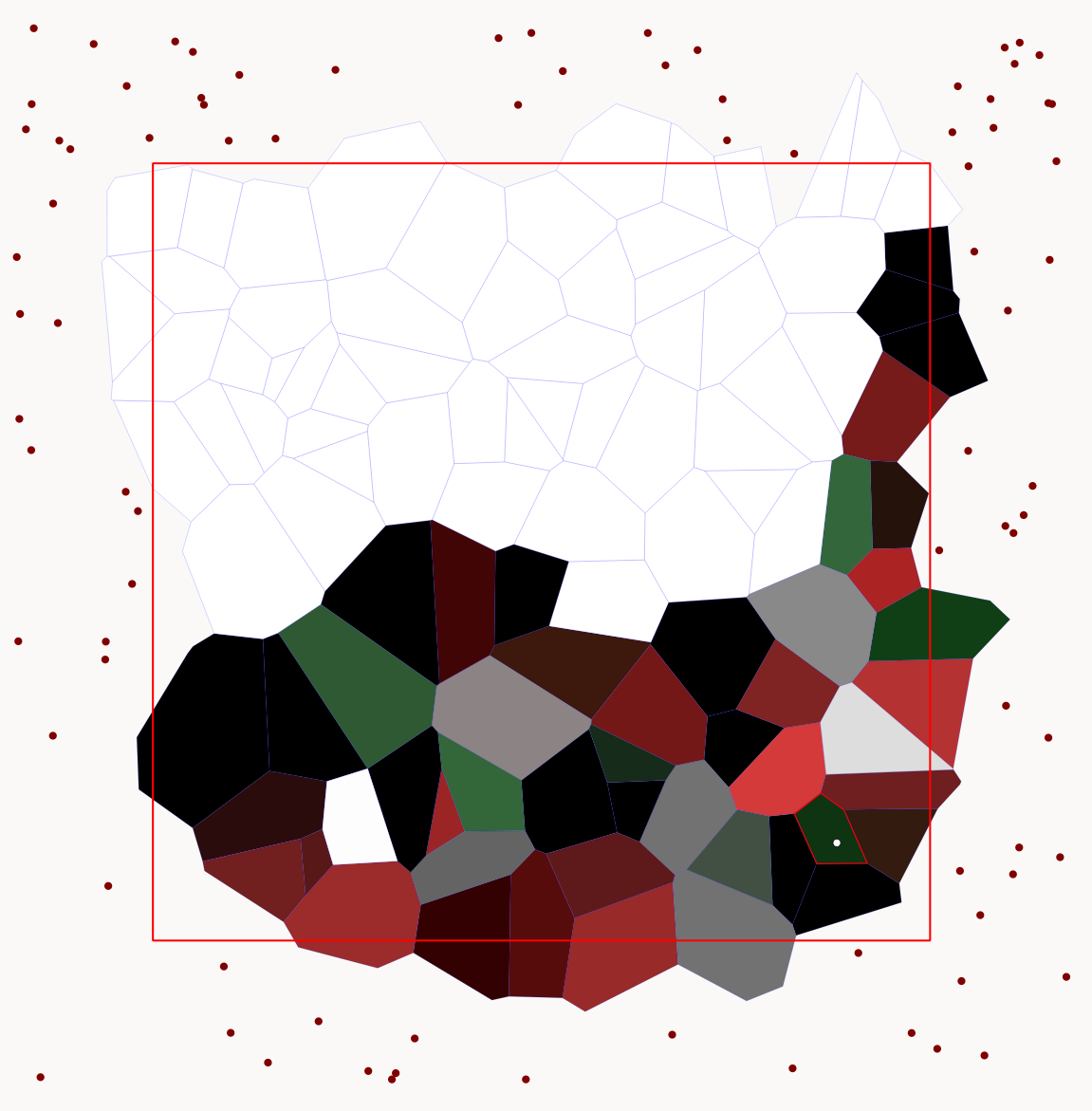}} \hspace*{2em}
    \subfloat[Near a diffuse wall\label{fig:detail3}]{\includegraphics[width=.4\textwidth]{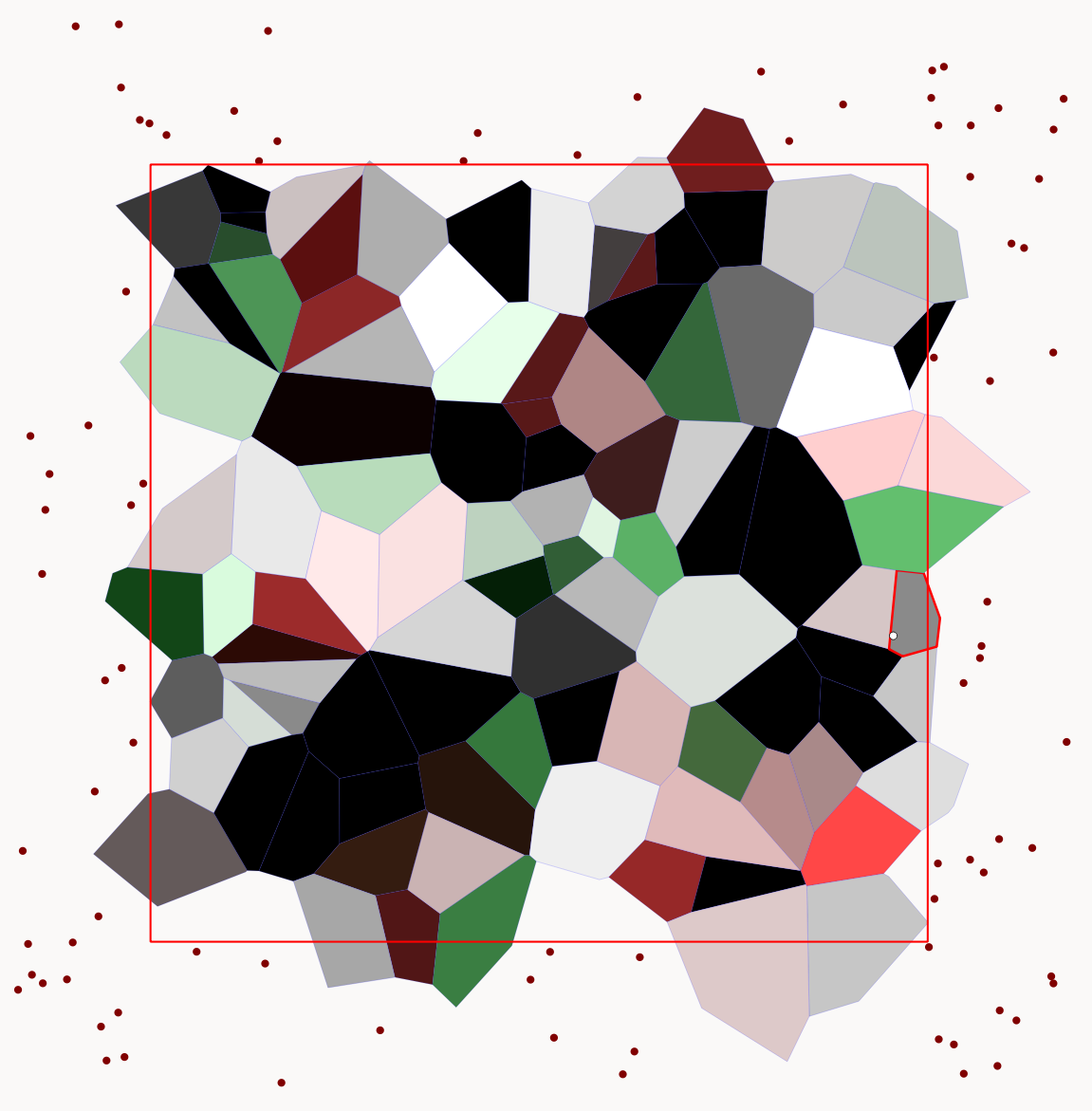}}%
    \caption{Pixels rendered with the $E_V$ estimators in a Cornell Box scene. Voronoi cells are colored according to the radiance of the sample computed by path tracing. Pixel boundaries are outlined in red. Voronoi weighting methods prove highly effective in scenarios where the radiance field exhibits significant spatial heterogeneity. They exhibit superior performance at major discontinuities, such as light source edges (figure a). However, they offer negligible improvement for diffuse materials (figure b), a limitation they share with adaptive sampling techniques.}
    \label{fig:trois_voronoi}
\end{figure}

Figure \ref{fig:main_figure} summarizes the key findings regarding the Voronoi estimators. The Voronoi estimators achieve superior performance by effectively reducing variance in regions exhibiting significant spatial radiance heterogeneity, specifically areas with steep radiance gradients. This improvement is particularly noticeable along light source boundaries, where the adaptive nature of the Voronoi weighting enables a more robust approximation of the radiance integral across sharp discontinuities like the edge of a projected light source. Additionally, this method also leads to faster convergence in areas involving specular light interactions on glossy materials. However, the benefits of the method are not uniform across the scene, especially where the radiance field is smooth or highly stochastic. Specifically, this method offers no substantial gain for surfaces with diffuse (Lambertian) reflection
or in regions exhibiting low or weak spatial correlation. This limitation is inherent to adaptive sampling techniques, which offer limited advantage over uniform sampling when the underlying light function is smooth or lacks structure \cite{Zwicker2015}. The geometry of the Voronoi tessellation for a direct light, a specular reflection or a diffuse wall is shown in Figure \ref{fig:trois_voronoi}.

\begin{figure}[htbp]
    \centering
    \subfloat[CPU time as a function of sample size\label{fig:graph1}]{ \includegraphics[width=.45\textwidth]{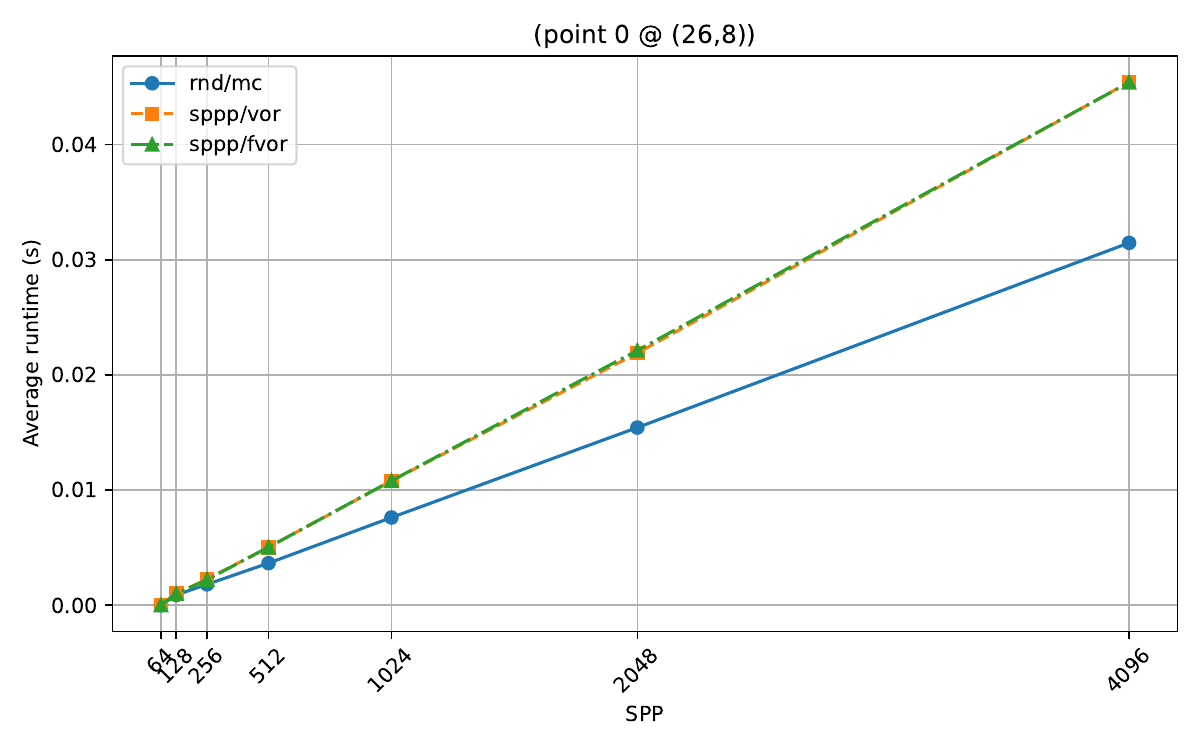}}
    \hfill
    \subfloat[variance as a function of CPU time\label{fig:graph2}]{\includegraphics[width=.45\textwidth]{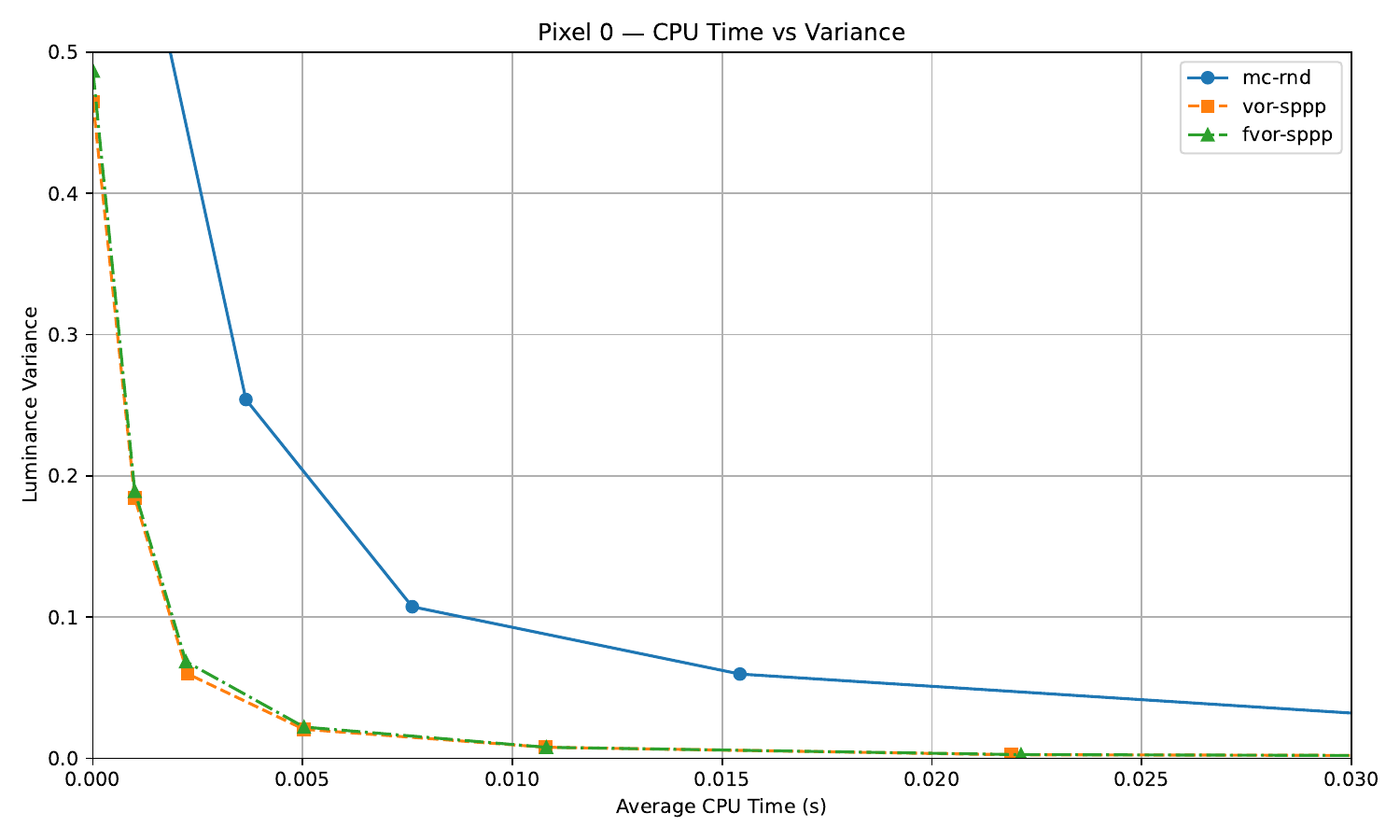}}%
    \caption{Estimators performance for a pixel located at the boundary of a light source. (a) Voronoi-based estimations lead to a nearly identical increase in the global computation time for both \texttt{vor} and \texttt{fvor} methods. (b) However, given a fixed CPU time budget, both methods prove to be more efficient.}
    \label{fig:runtime}
\end{figure}

Finally, the Voronoi-weighted integration method achieves significant reduction in variance in regions exhibiting steep radiance gradients, at the cost of some computational overhead. This overhead (see figure \ref{fig:graph1}) stems from constructing the Voronoi tessellation and computing the corresponding cell areas. Figure \ref{fig:graph2}, however, shows that this overhead is largely offset by the variance reduction achieved by the methods $E_V$ and $E_F$, making them more efficient than a standard Monte Carlo approach.

A more global performance indicator demonstrate an overall improvement of the numerical precision of the images computed by Voronoi estimators (Figure \ref{fig:main_figure}.) For a Cornell Box scene, Standard Monte Carlo estimation achieves an MSE value of 0.030, while a Voronoi estimation achieves a lower error (MSE=0.015).

\section{Discussion perspectives}

The integration techniques presented in this paper converge faster than a standard Monte Carlo approach for Hölder continuous functions. In practice, these methods prove more efficient than standard Monte Carlo when the integrand is computationally expensive.

Proposition \ref{prop:variance} is a new result regarding Poisson-Voronoi tessellations. In future work, this result could be extended to point processes that are not necessarily Poissonian and that include determinantal processes, in particular the Ginibre process. Using the notations of Proposition \ref{prop:variance}, a natural question that arises is whether 
$\sum_{x\in \eta_n}f(x)v_{\eta_n}(x)$ satisfies a central limit theorem, i.e. if 
\[ \textstyle \frac{\sum_{x\in \eta_n}f(x)v_{\eta_n}(x)-\EEE{\sum_{x\in \eta_n}f(x)v_{\eta_n}(x)}}{\VVV{\sum_{x\in \eta_n}f(x)v_{\eta_n}(x)}}\conv[n]{\infty}\mathcal{N}(0,1)\]
in distribution.

In the context of image rendering, our experiments highlighted that these techniques achieve rapid variance reduction in regions where the radiance field exhibits high-frequency structures. We demonstrated that a rejection scheme based on the dilation of the Voronoi seed space yields well-formed tessellations for this type of evaluation. Rather than relying on rejection, it would be of interest to investigate an adaptive approach; such a method would not reject an invalid tessellation but instead attempt to refine it by further expanding the integration window.

We restricted our analysis to window-like domains in dimension 2. Applying similar techniques beyond the primary ray would be a natural extension. Specifically, constructing a Voronoi tessellation within the directional hemisphere at each scene bounce may be worth considering. Furthermore, exploring Voronoi tessellations within the higher-dimensional space of paths represents another promising research perspective.

\medskip

\textbf{Contributions} All authors contributed to the study conception and design. Material preparation, data collection and analysis were performed by FV and SD. The first draft of the manuscript was written by FV and all authors commented on previous versions of the manuscript. All authors read and approved the final manuscript.

\bibliographystyle{plain}

\bibliography{biblio}

@book{Pharr2023,
  title={Physically based rendering: From theory to implementation},
  author={Pharr, Matt and Jakob, Wenzel and Humphreys, Greg},
  year={2023},
  month = {Mar},
  day = {28},
  address = "Cambridge, MA",
  publisher={MIT Press}
}

@BOOK{Robert2010,
  title     = "Monte Carlo Statistical Methods",
  author    = "Robert, Christian and Casella, George",
  publisher = "Springer",
  series    = "Springer texts in statistics",
  month     =  nov,
  year      =  2010,
  address   = "New York, NY",
  language  = "en"
}

@article{Zwicker2015,
author = {Zwicker, M. and Jarosz, W. and Lehtinen, J. and Moon, B. and Ramamoorthi, R. and Rousselle, F. and Sen, P. and Soler, C. and Yoon, S.-E.},
title = {Recent Advances in Adaptive Sampling and Reconstruction for Monte Carlo Rendering},
year = {2015},
issue_date = {May 2015},
publisher = {The Eurographs Association \& John Wiley \& Sons, Ltd.},
address = {Chichester, GBR},
volume = {34},
number = {2},
issn = {0167-7055},
journal = {Comput. Graph. Forum},
month = may,
pages = {667–681},
numpages = {15}
}

@book {Okabe2000,
    AUTHOR = {Okabe, A. and Boots, B. and Sugihara, K. and
              Chiu, S. N.},
     TITLE = {Spatial Tessellations: Concepts and Applications of {V}oronoi
              Diagrams},
    SERIES = {Wiley Series in Probability and Statistics},
   EDITION = {Second},
      PUBLISHER = {John Wiley \& Sons Ltd.},
   ADDRESS = {Chichester},
      YEAR = {2000},
     PAGES = {xvi+671},
      ISBN = {0-471-98635-6},
   MRCLASS = {52B55 (68U05)},
  MRNUMBER = {1770006 (2001c:52013)},
}

@Article{Kajiya1986,
  author    = {Kajiya, James T.},
  journal   = {ACM SIGGRAPH Computer Graphics},
  title     = {The rendering equation},
  year      = {1986},
  issn      = {0097-8930},
  month     = aug,
  number    = {4},
  pages     = {143--150},
  volume    = {20},
  doi       = {10.1145/15886.15902},
  publisher = {Association for Computing Machinery (ACM)},
}

@article {Calka02bis,
    AUTHOR = {Calka, P.},
     TITLE = {The distributions of the smallest disks containing the
              {P}oisson-{V}oronoi typical cell and the {C}rofton cell in the
              plane},
   JOURNAL = {Adv. in Appl. Probab.},
  FJOURNAL = {Advances in Applied Probability},
    VOLUME = {34},
      YEAR = {2002},
    NUMBER = {4},
     PAGES = {702--717},
      ISSN = {0001-8678},
     CODEN = {AAPBBD},
   MRCLASS = {60D05 (52C20 60G55)},
  MRNUMBER = {1938938 (2003j:60013)},
MRREVIEWER = {V. K. Oganyan},
       DOI = {10.1239/aap/1037990949},
       URL = {http://dx.doi.org/10.1239/aap/1037990949},
}

@article {Heveling09,
    AUTHOR = {Heveling, M. and Reitzner, M.},
     TITLE = {Poisson-{V}oronoi approximation},
   JOURNAL = {Ann. Appl. Probab.},
  FJOURNAL = {The Annals of Applied Probability},
    VOLUME = {19},
      YEAR = {2009},
    NUMBER = {2},
     PAGES = {719--736},
      ISSN = {1050-5164},
   MRCLASS = {60D05 (52A22 60C05 60G55 62F40)},
  MRNUMBER = {2521886 (2010d:60026)},
MRREVIEWER = {Lutz Muche},
       DOI = {10.1214/08-AAP561},
       URL = {http://dx.doi.org/10.1214/08-AAP561},
}

@article{Pages98,
title = {A space quantization method for numerical integration},
journal = {Journal of Computational and Applied Mathematics},
volume = {89},
number = {1},
pages = {1-38},
year = {1998},
issn = {0377-0427},
doi = {https://doi.org/10.1016/S0377-0427(97)00190-8},
url = {https://www.sciencedirect.com/science/article/pii/S0377042797001908},
author = {Gilles Pagès}
}

@Article{Reitzner12,
  author =       {Reitzner, M. and Spodarev, E. and Zaporozhetsz, D.},
  title =        {Set Reconstruction by {V}oronoi cells},
  journal =      { Adv. in Appl. Probab.},
  volume =       {44},
  pages =        {938--953},
  year =         {2012},
}

@Article{Reitzner24,
  author =       {Reitzner, M. and Strotmann, A.},
  title =        {Poisson-Delaunay approximation},
  journal =      { Available at https://arxiv.org/pdf/2410.23003.},
  year =         {2024},
}

@book {Schneider08,
    AUTHOR = {Schneider, R. and Weil, W.},
     TITLE = {Stochastic and Integral Geometry},
    SERIES = {Probability and its Applications (New York)},
 PUBLISHER = {Springer-Verlag},
   ADDRESS = {Berlin},
      YEAR = {2008},
     PAGES = {xii+693},
      ISBN = {978-3-540-78858-4},
   MRCLASS = {60-02 (52A22 60D05 60G55 62M30)},
  MRNUMBER = {2455326 (2010g:60002)},
MRREVIEWER = {V. K. Oganyan},
       DOI = {10.1007/978-3-540-78859-1},
}

@article {Schulte12bis,
    AUTHOR = {Schulte, M.},
     TITLE = {A central limit theorem for the {P}oisson-{V}oronoi
              approximation},
   JOURNAL = {Adv. in Appl. Math.},
  FJOURNAL = {Advances in Applied Mathematics},
    VOLUME = {49},
      YEAR = {2012},
    NUMBER = {3-5},
     PAGES = {285--306},
      ISSN = {0196-8858},
   MRCLASS = {60D05 (60F05 60G55)},
  MRNUMBER = {3017961},
       DOI = {10.1016/j.aam.2012.08.001},
       URL = {http://dx.doi.org/10.1016/j.aam.2012.08.001},
}

@article {Thale16,
    AUTHOR = {Th\"ale, C. and Yukich, J. E.},
     TITLE = {Asymptotic theory for statistics of the {P}oisson-{V}oronoi
              approximation},
   JOURNAL = {Bernoulli},
  FJOURNAL = {Bernoulli. Official Journal of the Bernoulli Society for
              Mathematical Statistics and Probability},
    VOLUME = {22},
      YEAR = {2016},
    NUMBER = {4},
     PAGES = {2372--2400},
      ISSN = {1350-7265,1573-9759},
   MRCLASS = {60G55 (60D05 60F25 62G20)},
  MRNUMBER = {3498032},
MRREVIEWER = {Markus\ Kiderlen},
       DOI = {10.3150/15-BEJ732},
       URL = {https://doi.org/10.3150/15-BEJ732},
}

@misc{Vandewiele2025,
  doi = {10.5281/zenodo.17910139},
  title     = "prise-3d/yapt: 1.0.3",
  author    = "Vandewiele, Franck and Delepoulle, Samuel",
  publisher = "Zenodo",
  year      =  2025
}

@inproceedings{Voechovsk2017,
  title = {Voronoi weighting of samples in Monte Carlo integration},
  DOI = {10.7712/120217.5385.17023},
  booktitle = {Proceedings of the 2nd International Conference on Uncertainty Quantification in Computational Sciences and Engineering (UNCECOMP 2017)},
  publisher = {Institute of Structural Analysis and Antiseismic Research School of Civil Engineering National Technical University of Athens (NTUA) Greece},
  address = {Rhodes Island, Greece},
  author = {Vořechovský, Miroslav and Sadílek, Václav and Eliáš, Jan},
  year = {2017},
  pages = {478--491},
}

@article{Guo2021,
  title = {Geometric Sample Reweighting for Monte Carlo Integration},
  volume = {40},
  ISSN = {1467-8659},
  url = {http://dx.doi.org/10.1111/cgf.14405},
  DOI = {10.1111/cgf.14405},
  number = {7},
  journal = {Computer Graphics Forum},
  publisher = {Wiley},
  author = {Guo,  J. and Eisemann,  E.},
  year = {2021},
  month = oct,
  pages = {109–119}
}

\appendix

\section{Proof of proposition \ref{prop:variance}}
\label{annex:proof}

\begin{proof}
Assertion \textit{(i)} directly follows from the Slivnyak-Mecke formula (see e.g. Corollary 3.2.3 in \cite{Schneider08}). Indeed, the latter implies
\begin{align*}
\EEE{\sum_{x\in \eta_n}f(x)v_{\eta_n}(x)} & = n\int_{\RR^d}f(x)\EEE{v_{\eta_n\cup\{x\}}(x)}\mathrm{d}x\\
& = \int_{\RR^d}f(x)\mathrm{d}x
\end{align*}
since $\EEE{v_{\eta_n\cup\{x\}}(x)} = \EEE{v_{\eta_n\cup\{0\}}(0)} = n^{-1}$.

 To deal with \textit{(ii)}, we apply the Poincaré inequality (see e.g. Section 2 in \cite{Reitzner12}), which claims that
\[\VVV{F(\eta_n)} \leq n \EEE{\int_{\RR^d} \left( F(\eta_n\cup\{x\}) - F(\eta_n)  \right)^2\mathrm{d}x  }\] for all measurable function $F:\mathbf{N}\rightarrow\RR_+$ such that $\EEE{F^2(\eta_n)}<\infty$, where $\mathbf{N}$ denotes the set of all locally finite counting measures in $\RR^d$. This implies

 \begin{equation}  \label{eq:poincarebound} \VVV{\sum_{x\in \eta_n}f(x)v_{\eta_n}(x)}
 \leq n\EEE{\int_{\RR^d} \left(  \sum_{y\in \eta_n\cup\{x\}} f(y)v_{\eta_n\cup\{x\}}(y) -  \sum_{y\in \eta_n} f(y)v_{\eta_n}(y) \right)^2   \mathrm{d}x}.\end{equation}



Now, let $x$ be fixed and let $y\in \eta_n$. Denote by $\mathcal{F}_{\eta_n}(y)$  the Voronoi flower associated with $y$, i.e.
\[\mathcal{F}_{\eta_n}(y) = \bigcup_{z\in C_{\eta_n}(y)}B(z,|z-y|),\] where $B(z,r)$ is the (closed) ball centered at $z$ with radius $r$.   Notice that, if $x\not\in \mathcal{F}_{\eta_n}(y)$, then $v_{\eta_n\cup\{x\}}(y)=v_{\eta_n}(y)$. Therefore

\begin{align*}
		\sum_{y\in \eta_n\cup\{x\}}  & f(y)v_{\eta_n\cup\{x\}}(y) - \sum_{y\in \eta_n} f(y)v_{\eta_n}(y) \\
	 &= \sum_{y\in \eta_n} f(y)\left( v_{\eta_n\cup\{x\}}(y) - v_{\eta_n}(y)\right) + f(x)v_{\eta_n\cup\{x\}}(x) \\
	 &= \sum_{y\in \eta_n} f(y)\left( v_{\eta_n\cup\{x\}}(y) - v_{\eta_n}(y)\right) \ind{x\in \mathcal{F}_{\eta_n}(y)} + f(x)v_{\eta_n\cup\{x\}}(x) \\
	 &= \sum_{y\in \mathcal{N}_{\eta_n\cup\{x\}}(x)} f(y)\left( v_{\eta_n\cup\{x\}}(y) - v_{\eta_n}(y)\right) + f(x)v_{\eta_n\cup\{x\}}(x),
	\end{align*}
where the last 
equality
comes from the fact that $x\in \mathcal{F}_{\eta_n}(y)$ iff $x\in \mathcal{N}_{\eta_n\cup\{x\}}(y)$, i.e. iff $y\in \mathcal{N}_{\eta_n\cup\{x\}}(x)$, and where $\mathcal{N}_{\eta_n\cup\{x\}}(x)$ denotes the set of neighbors of $x$ in the Voronoi tessellation associated with $\eta_n\cup\{x\}$. Furthermore, adapting an argument of Section 4 in \cite{Heveling09}, we have
\[v_{\eta_n\cup\{x\}}(x) = \sum_{y\in \mathcal{N}_{\eta_n\cup\{x\}}(x)}\left( v_{\eta_n}(y) - v_{\eta_n\cup\{x\}}(y)  \right).\]
Thus
\[\sum_{y\in \eta_n\cup\{x\}} f(y)v_{\eta_n\cup\{x\}}(y) -  \sum_{y\in \eta_n} f(y)v_{\eta_n}(y)  = \sum_{y\in \mathcal{N}_{\eta_n\cup\{x\}}(x)}(f(y)-f(x))\left( v_{\eta_n\cup\{x\}}(y) - v_{\eta_n}(y)  \right).\]
This together with \eqref{eq:poincarebound} gives
 \begin{align}
 \label{eq:boundvariance}
 \VVV{\sum_{x\in \eta_n}f(x)v_{\eta_n}(x)} & \leq n \int_{\RR^d} \EEE{\left(  \sum_{y\in \mathcal{N}_{\eta_n\cup\{x\}}(x)}   |f(y)-f(x)||  v_{\eta_n}(y) - v_{\eta_n\cup\{x\}}(y) |  \right)^2}\notag\\
 & =: r_1(n)+r_2(n)
 \end{align}
where 
\[r_1(n):=  n \int_{\RR^d} \EEE{ \sum_{y\in \mathcal{N}_{\eta_n\cup\{x\}}(x)}   |f(y)-f(x)|^2|  v_{\eta_n}(y) - v_{\eta_n\cup\{x\}}(y) |^2}\mathrm{d}x\]
and
 \begin{multline*}
 r_2(n):= n \int_{\RR^d} \mathbb{E}\Big[\sum_{(y_1,y_2)_{\neq }\in \mathcal{N}_{\eta_n\cup\{x\}}(x)^2}   |f(y_1)-f(x)||f(y_2)-f(x)|\\
 \times |  v_{\eta_n}(y_1) - v_{\eta_n\cup\{x\}}(y_1) | |  v_{\eta_n}(y_2) - v_{\eta_n\cup\{x\}}(y_2) |\Big]\mathrm{d}x.
 \end{multline*}


We provide below a bound for $r_1(n)$. According to the Slivnyak-Mecke formula, we have
 \[r_1(n) =  n^2 \int_{(\RR^d)^2} |f(y)-f(x)|^2\EEE{ |v_{\eta_n\cup\{y\}}(y) - v_{\eta_n\cup\{x,y\}}(y) |^2\ind{y\in \mathcal{N}_{\eta_n\cup\{x,y\}}(x)}  }\mathrm{d}y\mathrm{d}x\]
 Since $v_{\eta_n\cup\{x,y\}}(y)\leq v_{\eta_n}(y)$, we get
 \begin{align*}
  r_1(n) & \leq  4n^2 \int_{(\RR^d)^2} |f(y)-f(x)|^2 \EEE{v^2_{\eta_n\cup\{y\}}(y) \ind{y\in \mathcal{N}_{\eta_n\cup\{x,y\}}(x)}   }\mathrm{d}y\mathrm{d}x\\
  &  = 4 \int_{(\RR^d)^2} |f(y)-f(x)|^2 \EEE{v^2_{\eta\cup\{n^{1/d}y\}}(n^{1/d}y) \ind{n^{1/d}y\in \mathcal{N}_{\eta\cup\{n^{1/d}x,n^{1/d}y\}}(n^{1/d}x)}   }\mathrm{d}y\mathrm{d}x,
 \end{align*}
where the last equality comes from the fact that $\eta_n \overset{\mathcal{D}}{=}n^{-1/d}\eta$, where $\eta:=\eta_1$ is a PPP with intensity 1. Taking the change of variables $x'=n^{1/d}x$ and $y'=n^{1/d}y$ and applying the Hölder's inequality, we obtain
 \[
 r_1(n)  \leq  4n^{-2} \int_{(\RR^d)^2} |f(n^{-1/d}y')-f(n^{-1/d}x')|^2 \EEE{v^2_{\eta\cup\{y'\}}(y') \ind{y'\in \mathcal{N}_{\eta\cup\{x',y'\}}(x')}   }\mathrm{d}y'\mathrm{d}x',
 \]
and therefore
 \begin{multline*}
 r_1(n) \leq 4n^{-2} \int_{(\RR^d)^2} |f(n^{-1/d}y')-f(n^{-1/d}x')|^2 \EEE{v^4_{\eta\cup\{y'\}}(y')}^{1/2}\\\times  \PPP{y'\in \mathcal{N}_{\eta\cup\{x',y'\}}(x')}^{1/2}\mathrm{d}y'\mathrm{d}x'
 \end{multline*}

Since $\eta$ is stationary, we have $ \EEE{v^4_{\eta\cup\{y'\}}(y')} = \EEE{v^4(\mathcal{C})}<\infty$ for any $y'\in \RR^d$, where $\mathcal{C} \overset{\mathcal{D}}{=}C_{\eta\cup\{0\}}(0)$ and $v(\mathcal{C})$ denote the typical cell associated with $\eta$ and the volume of  $\mathcal{C}$, respectively. Writing $W_n:=n^{1/d}W$, we get
 \begin{align*}
 r_1(n) & \leq cn^{-2} \int_{(\RR^d)^2} |f(n^{-1/d}y')-f(n^{-1/d}x')|^2 \PPP{y'\in \mathcal{N}_{\eta\cup\{x',y'\}}(x')}^{1/2}\mathrm{d}y'\mathrm{d}x'\\
 & = 2c n^{-2} \int_{W_n}\left( \int_{\RR^d} |f(n^{-1/d}y')-f(n^{-1/d}x')|^2 \PPP{y'\in \mathcal{N}_{\eta\cup\{x',y'\}}(x')}^{1/2}\mathrm{d}y'\right)\mathrm{d}x',
 \end{align*}


where, in the last line, we used the fact that $f(n^{-1/d}y') = f(n^{-1/d}x')$ if $(x',y')\in (W_n^c)^2$ since $W$ is the support of $f$. The factor 2 comes from the fact that the integrand is symmetric w.r.t. $x'$ and $y'$ (since $y'\in \mathcal{N}_{\eta\cup\{x',y'\}}(x')$ iff $x'\in \mathcal{N}_{\eta\cup\{x',y'\}}(y')$).  Furthermore, for any locally finite subset $\chi\subset \RR^d$ in general position and for any $x\in \chi$, let $R_\chi(x)$ be the circumscribed radius associated with $C_\chi(x)$, i.e.
\[R_\chi(x) = \inf\{r>0: B(x,r) \supset C_\chi(x)\}.\]
Since $\mathcal{N}_\chi(x)\subset B(x,2R_\chi(x))$, we get
 \begin{multline*}
 r_1(n)\\
 \begin{split}& \leq c n^{-2} \int_{W_n}\left( \int_{\RR^d}|f(n^{-1/d}y')-f(n^{-1/d}x')|^2 \PPP{y'\in B(x',2R_{\eta\cup\{x',y'\}}(x'))}^{1/2}\mathrm{d}y'\right)\mathrm{d}x'\\
 & \leq cn^{-2-\tfrac{2\alpha}{d}}\int_{W_n}\left( \int_{\RR^d} |y'-x'|^{2\alpha}  \PPP{y'\in B(x',2R_{\eta\cup\{x'\}}(x'))}^{1/2}\mathrm{d}y'  \right)\mathrm{d}x',
 \end{split}
 \end{multline*}
where, in the last line, we used the facts that $f$ satisfies a Hölder condition with exponent $\alpha \in (0,1]$ and that $R_{\eta\cup\{x',y'\}}(x')\leq R_{\eta\cup\{x'\}}(x')$. Thanks to the stationarity of $\eta$, for any $x'\in \RR^d$, we have $R_{\eta\cup\{x'\}}(x') \overset{\mathcal{D}}{=}R(\mathcal{C})$, where $R(\mathcal{C})$ denotes the circumscribed radius of the typical cell. Therefore
\[ \int_{\RR^d} |y'-x'|^{2\alpha}  \PPP{y'\in B(x',2R_{\eta\cup\{x'\}}(x'))}^{1/2}\mathrm{d}y'  =  \int_{\RR^d}  |z|^{2\alpha} \PPP{R(\mathcal{C})\geq \frac{|z|}{2}}^{1/2}\mathrm{d}z. \]
It is known that the tail of $R(\mathcal{C})$ decreases exponentially fast to 0, see e.g. \cite{Calka02bis}, which proves that the last term is finite. Integrating over $x'\in W_n$, we get
\[r_1(n)\leq cn^{-1-\tfrac{2\alpha}{d}}.\]

We provide below an upper bound for $r_2(n)$. This will be sketched because we will proceed in the same spirit as we did for $r_1(n)$. First, according to the Slivnyak-Mecke formula, we write
\begin{align*}
r_2(n) & = n^3 \int_{(\RR^d)^3} |f(y_1)-f(x)| |f(y_2)-f(x)| \mathbb{E}\Bigg[ \Big| v_{\eta_n \cup\{x,y_1,y_2\}}(y_1) - v_{\eta_n \cup\{y_1,y_2\}}(y_1) \Big| \\
& \quad \times \Big| v_{\eta_n \cup\{x,y_1,y_2\}}(y_2) - v_{\eta_n \cup\{y_1,y_2\}}(y_2) \Big| 
\ind{(y_1,y_2)\in \mathcal{N}^2_{\eta_n\cup\{x,y_1,y_2\}}(x)} \Bigg] \mathrm{d}y_2\mathrm{d}y_1\mathrm{d}x.
\end{align*}
Taking the change of variables $x'=n^{1/d}x$, $y'_1=n^{1/d}y_1$, $y'_2=n^{1/d}y_2$ and using the facts that $\eta_n\overset{\mathcal{D}}{=}n^{-1/d}\eta$, $v_{\eta\cup\{x',y'_1,y'_2\}}(y'_i)\leq v_{\eta\cup\{y'_1,y'_2\}}(y'_i)$, with $i\in \{1,2\}$, we get
 \begin{multline*}
 r_2(n) \leq 4cn^3n^{-2}n^{-3}\int_{W_n} \Big( \int_{(\RR^d)^2} |f(n^{-1/d}y'_1)-f(n^{-1/d}x')||f(n^{-1/d}y'_2)-f(n^{-1/d}x')|\\
 \times   \EEE{v_{\eta\cup\{y'_1,y'_2\}}(y'_1)v_{\eta\cup\{y'_1,y'_2\}}(y'_2) \ind{(y'_1,y'_2)\in \mathcal{N}^2_{\eta\cup\{x',y'_1,y'_2\}}(x')} }\mathrm{d}y'_2\mathrm{d}y'_1\Big) \mathrm{d}x'.
 \end{multline*}
Notice that, without restriction, we have integrated over $x$ by using symmetry arguments and the fact that $f$ has compact support $W$. Furthemore, according to the Hölder's inequality, we have for any $(x,y'_1,y'_2)\in W_n\times \RR^d\times \RR^d$,
 \begin{multline*}
  \EEE{v_{\eta\cup\{y'_1,y'_2\}}(y'_1)v_{\eta\cup\{y'_1,y'_2\}}(y'_2) \ind{(y'_1,y'_2)\in \mathcal{N}^2_{\eta\cup\{x',y'_1,y'_2\}}(x')} }\\
  \begin{split} 
 & \leq \left(\EEE{v^4_{\eta\cup\{y'_1,y'_2\}}(y'_1)} \EEE{v^4_{\eta\cup\{y'_1,y'_2\}}(y'_1)} \PPP{(y'_1,y'_2)\in \mathcal{N}^2_{\eta\cup\{x',y'_1,y'_2\}}(x')}\right)^{1/4}\\
 & \leq c \left( \PPP{(y'_1,y'_2)\in \mathcal{N}^2_{\eta\cup\{x',y'_1,y'_2\}}(x')}\right)^{1/4},
 \end{split}
 \end{multline*}
where, in the last line, we used the facts that $v_{\eta\cup\{y'_1,y'_2\}}(y'_i)\leq v_{\eta\cup\{y'_i\}}(y'_i)\overset{\mathcal{D}}{=}v(\mathcal{C})$ and that $\EEE{v^4(\mathcal{C})}<\infty$. This, together with the facts that $f$ is Hölder and $\mathcal{N}^2_{\eta\cup\{x',y'_1,y'_2\}}(x')\subset B(x',2R_{\eta\cup\{x'\}}(x'))$ implies
 \begin{multline*}
 r_2(n)\leq cn^{-2-\tfrac{2\alpha}{d}}\int_{W_n} \Big( \int_{(\RR^d)^2}|x'-y'_1|^\alpha|x'-y'_2|^\alpha \\
 \times   \PPP{y'_1\in B(x',2R_{\eta\cup\{x'\}}(x'))}^{1/4}  \PPP{y'_2\in B(x',2R_{\eta\cup\{x'\}}(x'))}^{1/4} \mathrm{dy'_1}\mathrm{d}y'_2   \Big)\mathrm{d}x'.
 \end{multline*}
Following the same approach as for $r_1(n)$, we get
\begin{align*}
r_2(n) & \leq cn^{-2-\tfrac{2\alpha}{d}}  \int_{W_n}\left(  \int_{\RR^d} |z|^\alpha \PPP{R(\mathcal{C}) \geq \frac{|z|}{2}}^{1/4}  \right)^2\mathrm{d}x'\\
& \leq cn^{-1-\frac{2\alpha}{d}}.
\end{align*}
This concludes the proof of Proposition \ref{prop:variance}. 
\end{proof}

\end{document}